\documentclass[11pt]{amsart}
\usepackage{graphicx}
\usepackage{amssymb}
\usepackage{amsmath}
\usepackage{amsthm}
\usepackage{amsfonts}
\usepackage{bbm}
\usepackage{tikz}
\usepackage{tikz-cd}
\usepackage{setspace,kantlipsum}
\usepackage[toc,page]{appendix}
\usepackage{hyperref}
\usepackage{xcolor}
\usepackage{overarrows}
\usepackage{multirow}
\usepackage{bm}
\usepackage{caption}
\usepackage{placeins}

\definecolor{green}{RGB}{0, 204, 0}
\definecolor{forest}{RGB}{0, 118, 0}
\definecolor{orange}{RGB}{255, 170, 0}
\definecolor{blue}{RGB}{0, 64, 255}
\definecolor{purple}{RGB}{163, 0, 163}

\usepackage{float}

\newtheorem*{thm}{Theorem}

\newtheorem{thrm}{Theorem}[section]

\newtheorem{lemma}[thrm]{Lemma}

\newtheorem*{quest}{Question}

\theoremstyle{definition}

\newtheorem{defn}[thrm]{Definition}

\newcommand{\id}{\operatorname{id}}

\makeatletter
\newcommand{\colim@}[2]{%
  \vtop{\m@th\ialign{##\cr
    \hfil$#1\operator@font colim$\hfil\cr
    \noalign{\nointerlineskip\kern1.5\ex@}#2\cr
    \noalign{\nointerlineskip\kern-\ex@}\cr}}%
}
\newcommand{\colim}{%
  \mathop{\mathpalette\colim@{\rightarrowfill@\scriptscriptstyle}}\nmlimits@
}
\renewcommand{\varprojlim}{%
  \mathop{\mathpalette\varlim@{\leftarrowfill@\scriptscriptstyle}}\nmlimits@
}
\renewcommand{\varinjlim}{%
  \mathop{\mathpalette\varlim@{\rightarrowfill@\scriptscriptstyle}}\nmlimits@
}
\makeatother

\usepackage[margin=1in]{geometry}

\title{Infinite order corks that admit Stein structures}
\author[Y. Teng]{Yikai Teng} \address{Rutgers University-Newark,
Newark, NJ 07103} \email{yikai.teng@rutgers.edu}
\begin{document}
\begin{spacing}{1}

\begin{abstract}
    In this paper, we construct an infinite family of infinite order corks admitting Stein structures. This answers a question raised by Akbulut \cite{Akb16b,Akb24} and Gompf \cite{Gom17a}.\vspace{-25pt}
\end{abstract}

\maketitle

\section{Introduction}

    A cork is a pair $(W,f)$, where $W$ is a compact, contractible, smooth 4-manifold, and $f:\partial W\rightarrow \partial W$ is a boundary diffeomorphism that does not extend to a diffeomorphism $W\rightarrow W$. In \cite{CFHS} and \cite{Mat96}, it is proved that any exotic copy of a simply-connected closed 4-manifold can be obtained by a single cut-and-paste along an embedded cork. Moreover, in \cite{AM98}, the above cork twisting theorem was upgraded such that the cork chosen can be modified to be a Stein domain. 

    In \cite{Gom17b} and \cite{Gom17a}, Gompf constructed a family of examples of infinite order corks, namely a pair $(W,f)$ such that for any positive integer $k$, the pair $(W, f^k)$ defines a non-trivial cork. To be specific and for future convenience, we give the following definition.

    \begin{defn}
        A smooth 4-manifold $X$ with nonempty boundary is said to \textit{admit an infinite order boundary twist} if there exists a boundary automorphism $\phi: \partial X \rightarrow \partial X$ such that the iterated automorphism $\phi^k$ does not extend to a self-diffeomorphism of the 4-manifold $X$ for any $k>0$.
    \end{defn}

    However, it is not known whether there exist examples of infinite order corks that admit Stein structures. Gompf, in \cite{Gom17a}, and Akbulut, in \cite{Akb16b} and \cite{Akb24}, asked the question whether any of Gompf's examples can be equipped with a Stein structure. Our main result gives a positive answer to this question.

    \begin{thm}{\label{thm-cork-stein}}
        Gompf's cork $C(1,1;-1)$ is a Stein domain. Moreover, the cork $C(1,1;-1)$ can be modified into an infinite family of infinite order corks that all admit Stein structures. 
    \end{thm}

\begin{figure}[H]
                        \centering
                        \resizebox{11.7cm}{!}{
                        \begin{picture}(15.0cm,5.4cm)
                            \put(0cm,0.0cm){\includegraphics[width=6cm, height=5.4cm]{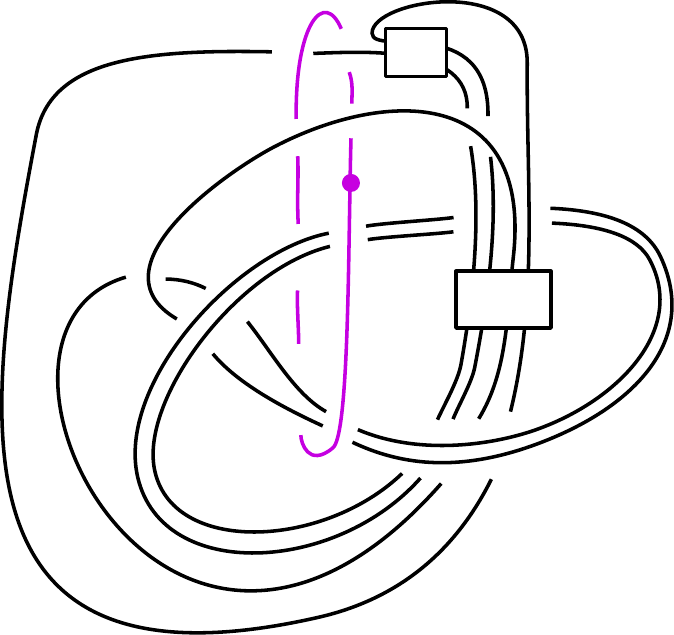}}
                            \put(8.8cm,0cm){\includegraphics[width=6.2cm, height=5.4cm]{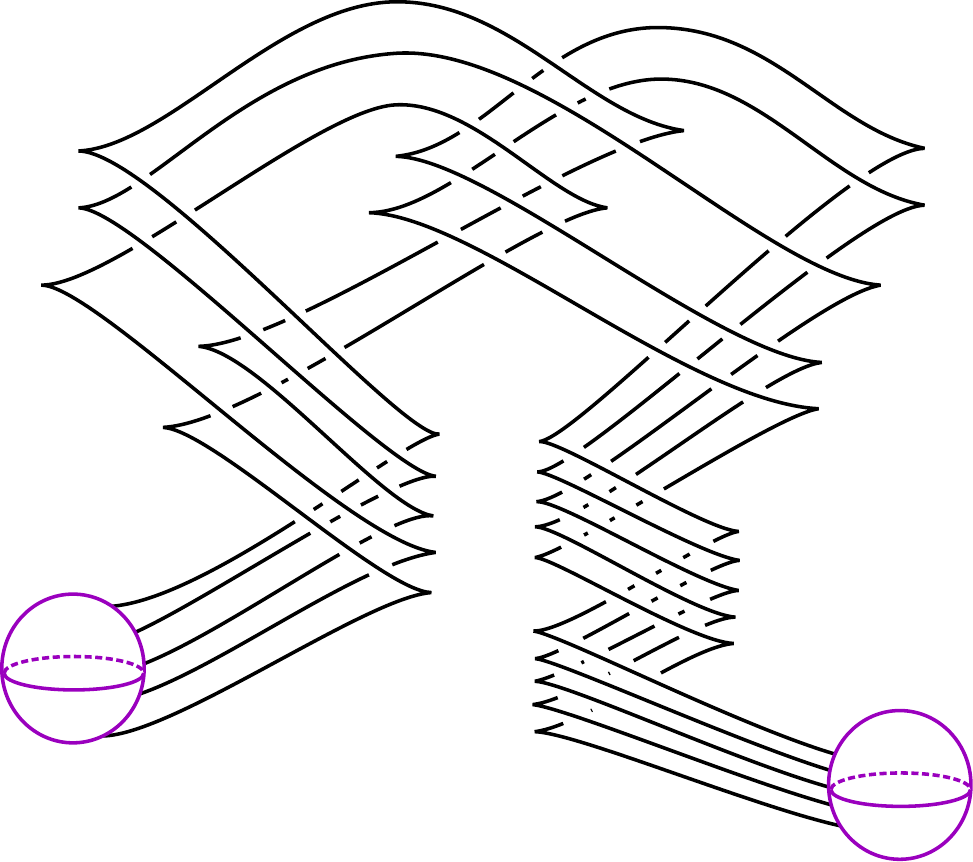}}
                            \put(3.47cm,4.85cm){\footnotesize{$-3$}}
                            \put(4.22cm,2.74cm){$-1$}

                            \put(4.85cm,4.32cm){$0$}
                            \put(8.95cm,4.82cm){$-2$}
                        
                        \end{picture} 
                        }
                        \caption{Left: Gompf's cork $C(1,1;-1)$ interpreted as a Mazur-type manifold. Right: a Legendrian diagram of $C(1,1;-1)$ whose 2-handle has framing $-2$ and Thurston-Bennequin number $-1$. }
                        \label{20240604-1}
            \end{figure}
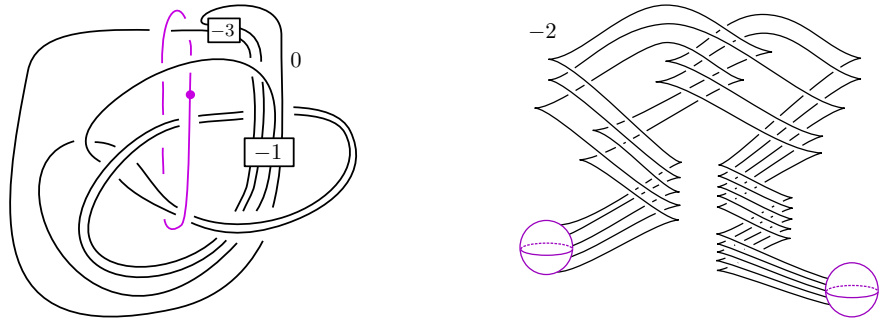

    To prove the first assertion, we will begin by realizing $C(1,1;-1)$ as a Mazur-type manifold, that is, as a manifold that admits a handle decomposition with a single $1$-handle and a single $0$-framed $2$-handle. The Kirby diagram is shown on the left in Figure \ref{20240604-1}. We will then modify this diagram into the Legendrian handle diagram shown on the right in Figure \ref{20240604-1}. By Gompf's handlebody criterion for Stein surfaces \cite{Gom98}, this Legendrian handle diagram defines a Stein structure on $C(1,1;-1)$. 

    To construct an infinite family of infinite order Stein corks from this single example, we will prove a diagrammatic auxiliary lemma: if the dotted circle representing the $1$-handle in Figure \ref{20240604-1} is replaced by any slice knot contained in its tubular neighborhood, the resulting 4-manifold also admits an infinite order boundary twist. The details will be given in Section \ref{nutshell}.

\begin{figure}[H]
                        \centering
                        \resizebox{4.95cm}{!}{
                        \begin{picture}(5.5cm,4cm)
                            \put(0cm,0cm){\includegraphics[width=5.5cm, height=4cm]{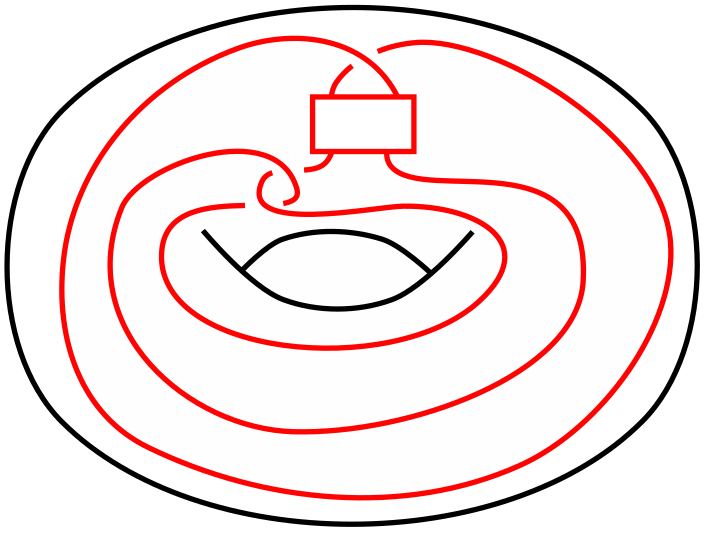}}
                            \put(2.5cm,3.0cm){\textcolor{red}{$-m$}}
                        \end{picture}
                        }
                        \caption{Yasui's family of patterns $P_m$, which will be used to modify $C(1,1;-1)$.}
                        \label{20240604-2}
                \end{figure}

\begin{figure}[H]
                        \centering
                        \resizebox{10.8cm}{!}{
                        \begin{picture}(12cm,9cm)
                            \put(0cm,0cm){\includegraphics[width=12cm, height=9cm]{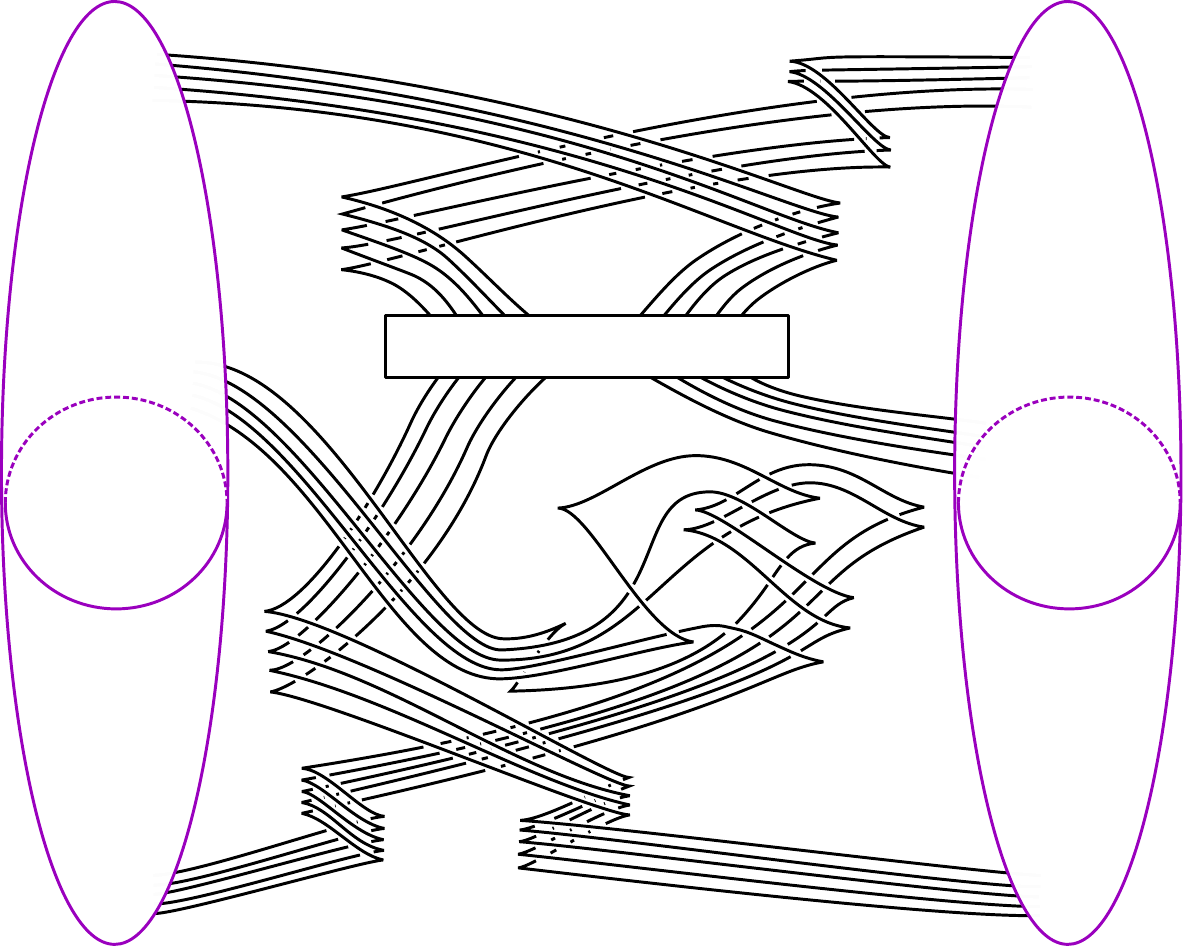}}
                            
                            \put(5.62cm,5.57cm){$-m$}
                            \put(3.45cm,8.40cm){$0$}
                        
                        \end{picture}    
                        }
                        \caption{Legendrian diagrams for the family of corks $C_m$, whose 2-handle has framing $0$ and Thurston-Bennequin number $+1$. }
                        \label{20240605-0}
            \end{figure}

    Using this auxiliary lemma, we will construct an infinite family of infinite order corks $C_m$ by replacing the purple 1-handle curve with its satellite using Yasui's generalized Mazur pattern $P_m$ \cite{Yas17}, as shown in Figure \ref{20240604-2}. In fact, each $C_m$ admits a Stein structure. Their Stein handle diagrams are shown in Figure \ref{20240605-0}.

    On the other hand, it remains an open question whether the boundary diffeomorphisms of these Stein corks are contactomorphisms, hence making them ``contact corks''.

    \begin{quest}
        For the Stein corks $C(1,1;-1)$ and $C_m$, are the boundary diffeomorphisms described in \cite{Gom17b} and Section \ref{nutshell} contactomorphisms?
    \end{quest}

    \subsection{Organization of the paper.}

        In Section \ref{nutshell}, we recall Gompf's construction of $C(1,1;-1)$ from \cite{Gom17b}, \cite{Gom17a} and describe the modifications used to construct $C_m$, with an emphasis on their Kirby diagrams. Using these diagrams, we prove in Section \ref{section-original-stein} that Gompf's cork $C(1,1;-1)$ admits a Stein structure, and in Section \ref{section-modified-stein} that the modified corks $C_m$ admit Stein structures. Together, these results prove the main theorem.  

    \subsection{Acknowledgments.}

        The author is grateful to Kyle Hayden and Arunima Ray for suggesting this problem and offering continued guidance.

\section{Gompf's construction in a nutshell}{\label{nutshell}}

    \subsection{The essence of Gompf's construction}{\label{section-understand-gompf}}

\begin{figure}[H]
            \centering
            \resizebox{7.2cm}{!}{
            \begin{picture}(8cm,6cm)
                \put(0cm,0cm){\includegraphics[width=8cm, height=6cm]{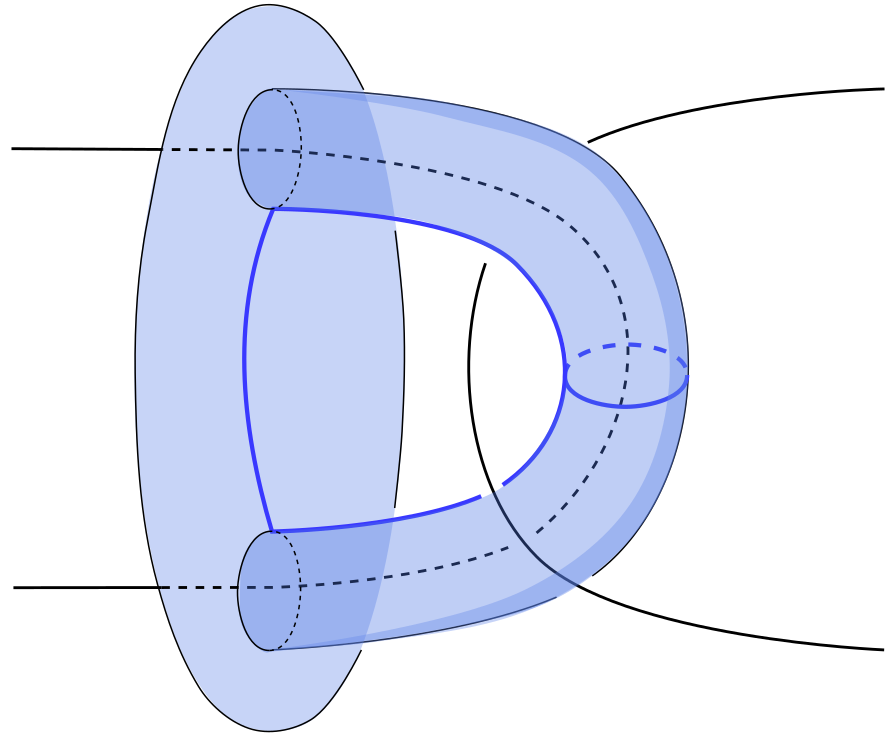}}
                \put(6.9cm,0.9cm){$K_k$}
                \put(1.0cm,5.2cm){$\partial \Sigma$}
                \put(4.6cm,0.65cm){$\Sigma$}
                \put(2.4cm,2.9cm){\textcolor{blue}{$C_{+1}$}}
                \put(6.3cm,2.9cm){\textcolor{blue}{$C_{-1}$}}
                        
            \end{picture}
            }
            \caption{The local picture for $M_k$, the complement of the $k$-twist knot, near its clasp.}
            \label{inf-clasp}
        \end{figure}
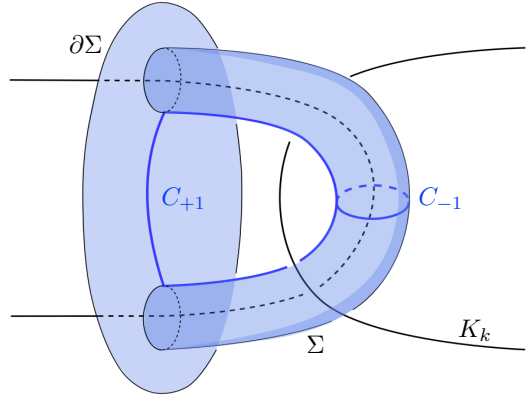

    In this subsection, we review the core idea of Gompf's construction. Figure \ref{inf-clasp} depicts a local picture of the $k$-twist knot (recall that the 0-twist knot is just the unknot drawn in a nonstandard diagram) near its unique clasp. We denote its complement by $M_k$.

\begin{figure}[H]
                        \centering
                        \resizebox{3.6cm}{!}{
                        \begin{picture}(4cm,4cm)
                            \put(0cm,0cm){\includegraphics[width=4cm, height=4cm]{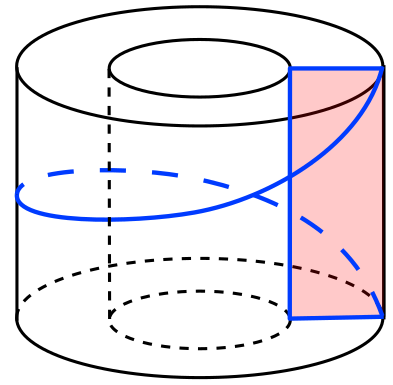}}
                            \put(3.6cm,0cm){$\partial \Sigma$}
                            \put(-0.2cm,1.95cm){$I$}
                        
                        \end{picture}   
                        }
                        \caption{The schematic picture for the Dehn twist $g$.}
                        \label{dehn-g}
            \end{figure}

    Let $\Sigma$ be the punctured torus shown in Figure \ref{inf-clasp}. Recall that the $k$-twist knot complement $M_k$ can be obtained by performing a $-\frac{1}{k}$-Dehn surgery along $\partial \Sigma$. Equivalently, $M_k$ can be obtained from $M_0$ by slitting $M_0$ open along $I\times\partial\Sigma$ and regluing by $g^k$, where $g$ is the Dehn twist shown in Figure \ref{dehn-g}. In this figure, the horizontal annulus is a collar neighborhood of $\partial \Sigma$ in $\Sigma$, and the meridian of the solid torus, which bounds the red disk in the picture, is sent by the Dehn twist $g$ to the blue curve. 

        Next, we pass to dimension 4 by taking the product with $S^1$. After identifying $M_0\times S^1$ with a neighborhood of a regular fiber in the elliptic surface $E(n)$ where $n\geq 2$, the process of cutting $E(n)$ open along $N:=I\times \partial\Sigma\times S^1$ and regluing by $g^k\times \id_{S^1}$ is equivalent to replacing the chosen fiber neighborhood with $M_k\times S^1$, that is, a Fintushel-Stern knot surgery \cite{FS98} on $E(n)$ by the $k$-twist knot.

        Now let $W$ be any embedded codimension-$0$ submanifold of $E(n)$ such that $N$ is contained in $\partial W$. We can then define a boundary diffeomorphism of $W$ by extending $g\times \id_{S^1}$ by the identity over the rest of $\partial W$. By the preceding analysis, cutting out the embedded 4-manifold $W$ from $E(n)$ and gluing it back in via the $k$-fold iterate of this boundary diffeomorphism yields the manifold obtained from $E(n)$ by Fintushel-Stern knot surgery using the $k$-twist knot. Since the $k$-twist knots are distinguished by their Alexander polynomials, the resulting knot-surgered elliptic surfaces are pairwise nondiffeomorphic \cite{FS98}. It follows that the embedded submanifold $W$ admits an infinite order boundary twist.

        One such example of $W$ is the $4$-manifold $I\times \Sigma\times S^1 \subset E(n)$. This infinite order manifold is the core of Gompf's construction. In particular, its boundary decomposes into two parts: 
        \begin{enumerate}
            \item the \textit{essential region}, namely $N = I\times \partial \Sigma \times S^1$, which contains the support of the infinite order boundary twist.
            \item the \textit{``playground''}, namely $\partial I\times \Sigma\times S^1$. 
        \end{enumerate}

        More importantly, we can modify the 4-manifold $I\times \Sigma\times S^1$ freely (attach handles, carve out surfaces, etc.) without affecting its infinite order boundary twist property, provided that 
        \begin{enumerate}
            \item any modifications on the boundary are supported in the playground region; and
            \item the resulting 4-manifold still admits an embedding in $E(n)$ whose restriction to a neighborhood of $N$ agrees with the original embedding.
        \end{enumerate}

    \subsection{Kirby diagrams}

    In this subsection, we construct Gompf's infinite order cork $C(1,1;-1)$ starting from the core manifold $I\times \Sigma \times S^1$, primarily in terms of Kirby diagrams, to establish the starting point for the calculations in Sections \ref{section-original-stein} and \ref{section-modified-stein}.

\begin{figure}[H]
            \centering
            \resizebox{13.6cm}{!}{
                \begin{picture}(17cm,5cm)
                    \put(0cm,0cm){\includegraphics[width=5cm, height=5cm]{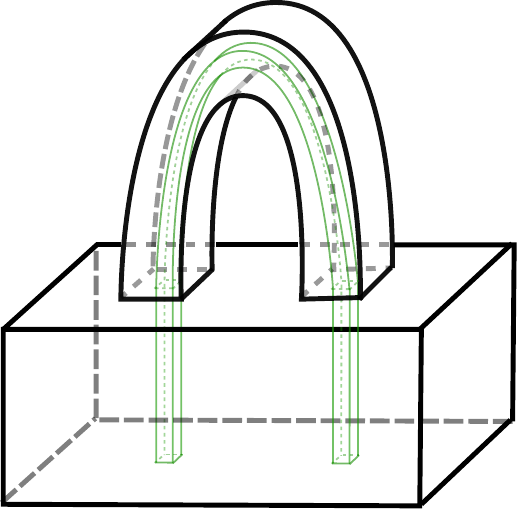}}
        
                    \put(6cm,0cm){\includegraphics[width=5cm, height=5cm]{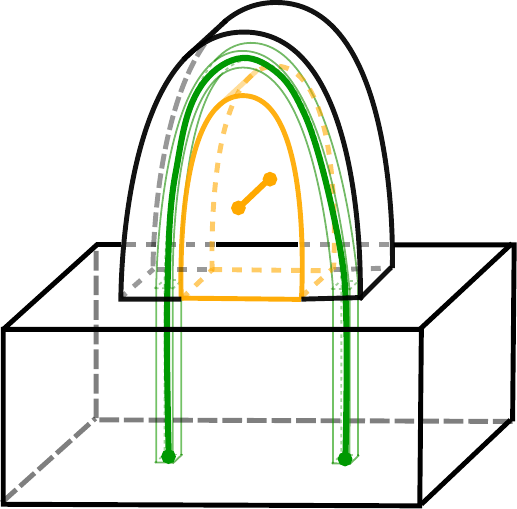}}
        
                    \put(12cm,0cm){\includegraphics[width=5cm, height=5cm]{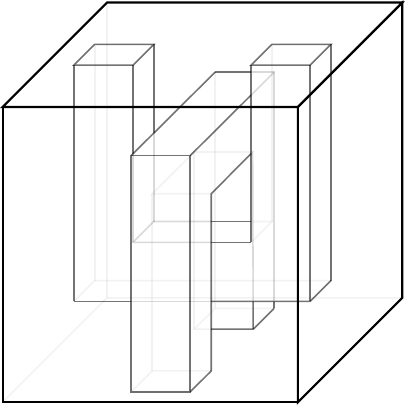}}
                                
                \end{picture}
            }
            \caption{Left: the $3$-manifold $\Sigma\times I$ may be viewed as a $3$-ball with a single $3$-dimensional $1$-handle attached, with the core of the handle removed and extended through $B^3$. Middle: after introducing a canceling $1/2$-handle pair, the same $3$-manifold can be viewed as the complement of the tubular neighborhoods of the green and yellow curves in $B^3$. Right: after an isotopy, this becomes the standard picture of a once-clasped $4$-ended tangle in $B^3$.}
            \label{clasp-complement}
        \end{figure}

\begin{figure}
        \centering
        \resizebox{12.15cm}{!}{
        \begin{picture}(13.5cm,4cm)
            \put(0cm,0cm){\includegraphics[width=6cm, height=4cm]{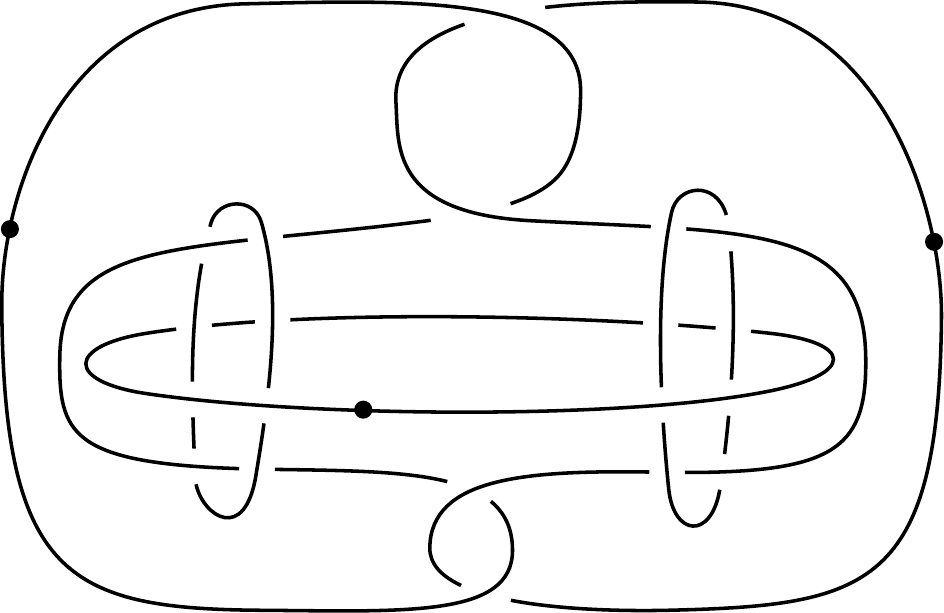}}
            \put(4.6cm, 2.85cm){$0$}
            \put(1.2cm, 2.8cm){$0$}

            \put(7.5cm,0cm){\includegraphics[width=6cm, height=4cm]{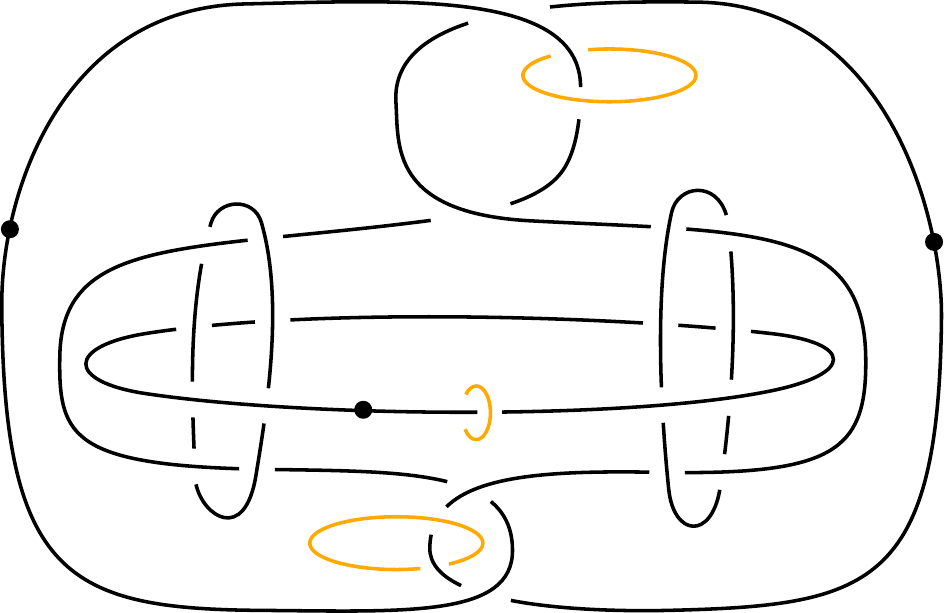}}
            \put(12.1cm, 2.85cm){$0$}
            \put(8.7cm, 2.8cm){$0$}
            \put(8.9cm, 0.25cm){\textcolor{orange}{$-1$}}
            \put(10.2cm, 1.6cm){\textcolor{orange}{$-1$}}
            \put(12.0cm, 3.4cm){\textcolor{orange}{$-1$}}

        \end{picture}
        }
        \caption{Left: The 4-manifold $I\times \Sigma\times S^1$. Right: The three $2$-handles are attached along the yellow curves.}
        \label{20241027-1}
    \end{figure}
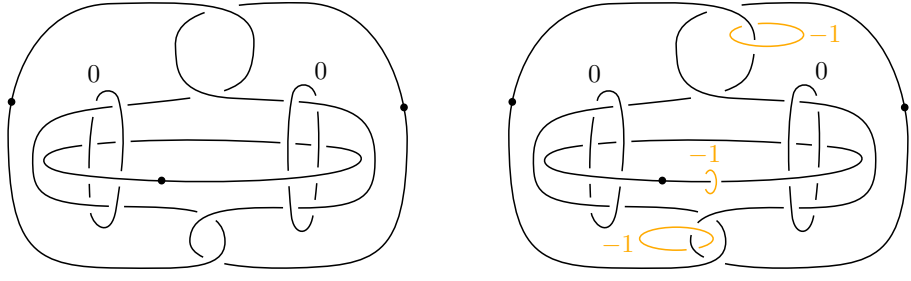

    We begin by considering the 3-manifold $I\times \Sigma$, where $\Sigma$ denotes a punctured torus. This 3-manifold is actually diffeomorphic to the clasp complement in $B^3$ shown on the right of Figure \ref{clasp-complement}. The rest of Figure \ref{clasp-complement} illustrates this identification. Following the standard procedure of drawing a Kirby diagram of (3-manifold)$\times S^1$ (cf. \cite[Section 6.5]{Akb16}), we obtain a handle diagram for the 4-manifold $I\times \Sigma \times S^1$, as shown on the left in Figure \ref{20241027-1}.

    Next, we modify the 4-manifold $I\times \Sigma\times S^1$ to obtain Gompf's cork $C(1,1;-1)$. We first attach $(-1)$-framed 2-handles along the curves $\{\pm 1\}\times C_{\pm 1} \times \{\theta_{\pm 1}\}$ and $\{+1\}\times \{p\} \times S^1$, where $\theta_{\pm 1}$ are two distinct generic points on $S^1$ and $p$ is a generic point in the interior of $\Sigma$. This gives the Kirby diagram on the right in Figure \ref{20241027-1}. 
    
    To see that this modified manifold still admits an infinite order boundary twist, we first observe that the attaching curves of the three $2$-handles lie entirely in the playground region of the boundary, and hence are disjoint from the essential region. On the other hand, recall that in the definition of knot surgery, the meridian of the surgery knot and the $S^1$-factor are identified with the two circle factors of a regular fiber of the elliptic fibration. Consequently, each of these circles comes with $6n$ parallel copies of the vanishing cycles in $E(n)$. In the present modification of $I\times \Sigma\times S^1$, the three $(-1)$-framed $2$-handles are attached precisely along these vanishing cycles, so the modified manifold still embeds in $E(n)$. Hence, the criteria from the previous subsection imply that this modified manifold still admits an infinite order boundary twist.
    
    We then drill out the two core disks of the 2-handles attached along $\{\pm 1\}\times C_{\pm 1} \times \{\theta_{\pm 1}\}$, extending them to the boundary by the annuli $I\times C_{\pm 1}\times \{\theta_{\pm 1}\}$. As before, the induced boundary modification is disjoint from the essential region. Moreover, since this operation merely carves out properly embedded disks, the resulting $4$-manifold still embeds in the elliptic surface $E(n)$. Hence, by the criteria from the previous subsection, this modified manifold still admits an infinite order boundary twist. This is Gompf's infinite order cork $C(1,1;-1)$.
    
\begin{figure}[H]
        \centering
        \resizebox{9.9cm}{!}{
        \begin{picture}(11cm,4.5cm)
            \put(0cm,0cm){\includegraphics[width=11cm]{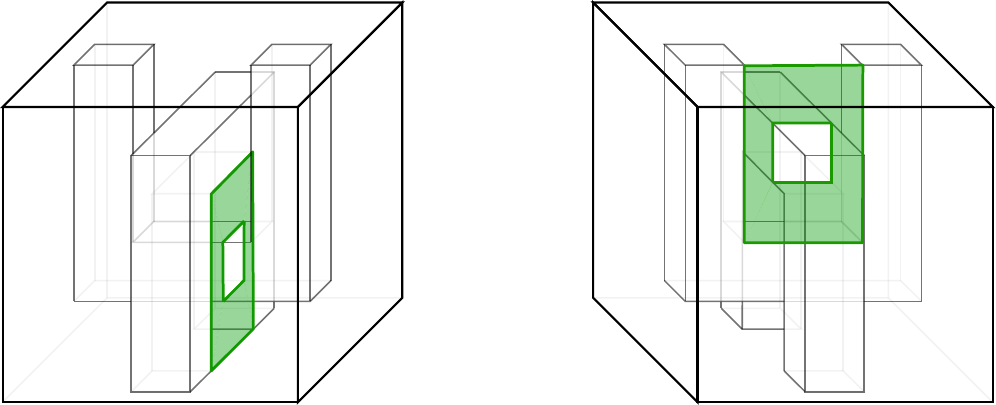}}
        \end{picture}
        }
        \caption{Deleting the cores of the 2-handles, extended by the green annuli, gives the infinite order cork $C(1,1;-1).$}
        \label{movie-with-bad-disks}
    \end{figure}

\begin{figure}[H]
        \centering
        \resizebox{12.15cm}{!}{
        \begin{picture}(13.5cm,4cm)
            \put(0cm,0cm){\includegraphics[width=6cm, height=4cm]{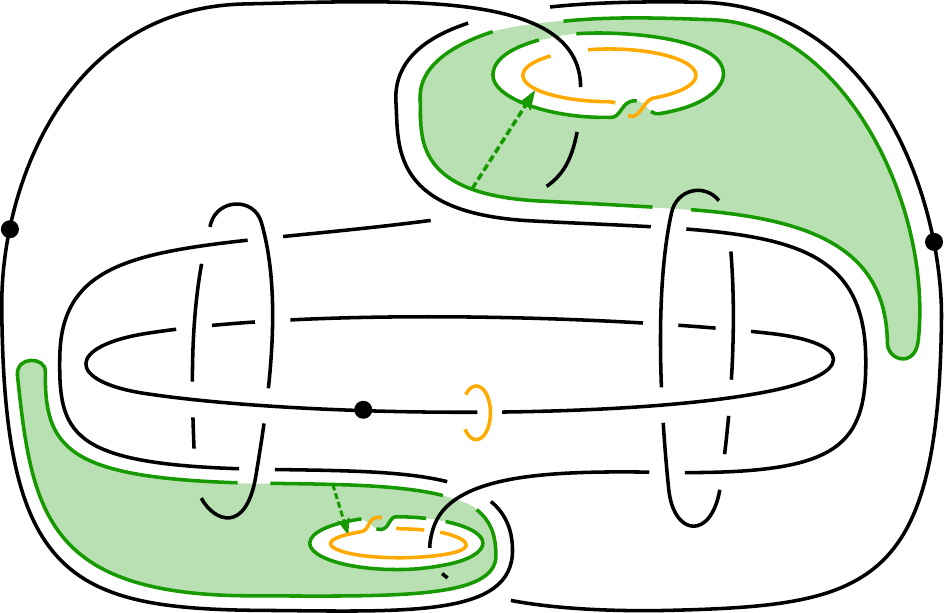}}
            \put(7.5cm,0cm){\includegraphics[width=6cm, height=4cm]{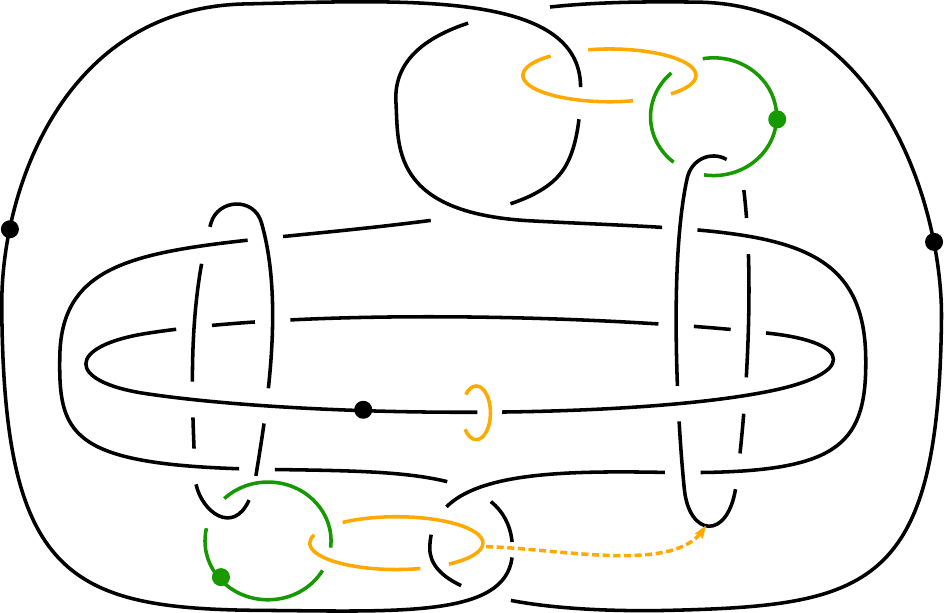}}
            \put(4.6cm, 2.85cm){$0$}
            \put(1.2cm, 2.8cm){$0$}
            \put(1.4cm, 0.25cm){\textcolor{orange}{$-1$}}
            \put(2.7cm, 1.6cm){\textcolor{orange}{$-1$}}
            \put(4.62cm, 3.4cm){\textcolor{orange}{$-1$}}

            \put(11.9cm, 3.07cm){$0$}
            \put(8.7cm, 2.8cm){$0$}
            \put(8.9cm, 0.35cm){\textcolor{orange}{$-1$}}
            \put(10.2cm, 1.6cm){\textcolor{orange}{$-1$}}
            \put(10.26cm, 3.4cm){\textcolor{orange}{$-1$}}

        \end{picture}   
        }
        \caption{Left: the disks to be drilled out, which consist of the green annuli, capped off by the cores of the yellow 2-handles. Right: the complement of the disks.\vspace{-7pt}}
        \label{20241027-2}
    \end{figure}
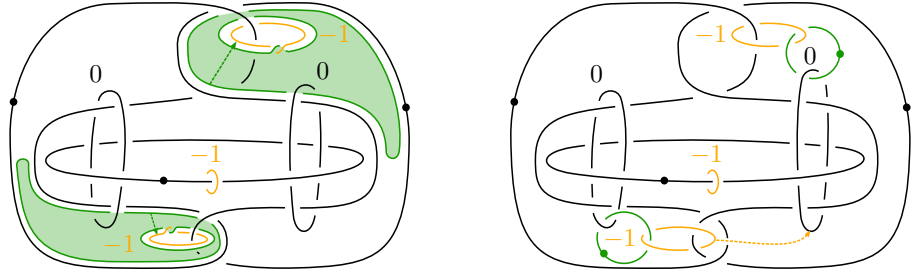

    Next, we describe this modification in Kirby-diagram language. We begin by locating the disks to be drilled out. Figure \ref{movie-with-bad-disks} records the curves $ \{\pm 1\}\times C_{\pm 1}\times \{\theta_{\pm1}\}$ together with the annuli they cobound. Translating to Kirby diagrams, we obtain the first picture of Figure \ref{20241027-2}, where the drilled-out surfaces are represented by the green annuli, capped off by the cores of the corresponding 2-handles. Indeed, the boundary curves of the annuli in Figure \ref{movie-with-bad-disks} correspond to the green curves in the first diagram in Figure \ref{20241027-2}: two are parallel copies of the attaching curves of the corresponding $2$-handles (so they are capped off by the core disks), while the other two run along longitudes of the opposite clasps. Following the procedure of \cite[Section 6.2]{GS99} (in particular, Figure 6.37), we obtain a Kirby diagram of the complement of these disks by sliding each disk off the corresponding $2$-handle along the indicated green arrows. After an isotopy, this gives the Kirby diagram on the right in Figure \ref{20241027-2}. 

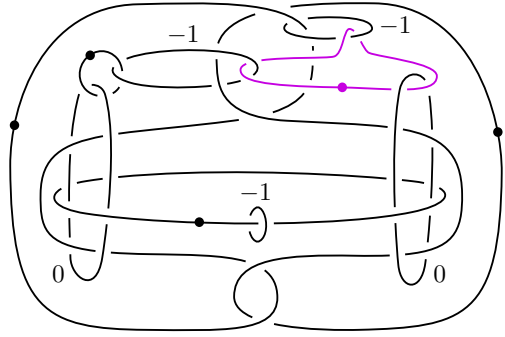
\begin{figure}[H]
        \centering
        \resizebox{6.75cm}{!}{
        \begin{picture}(7.5cm,5cm)
            \put(0cm,0cm){\includegraphics[width=7.5cm, height=5cm]{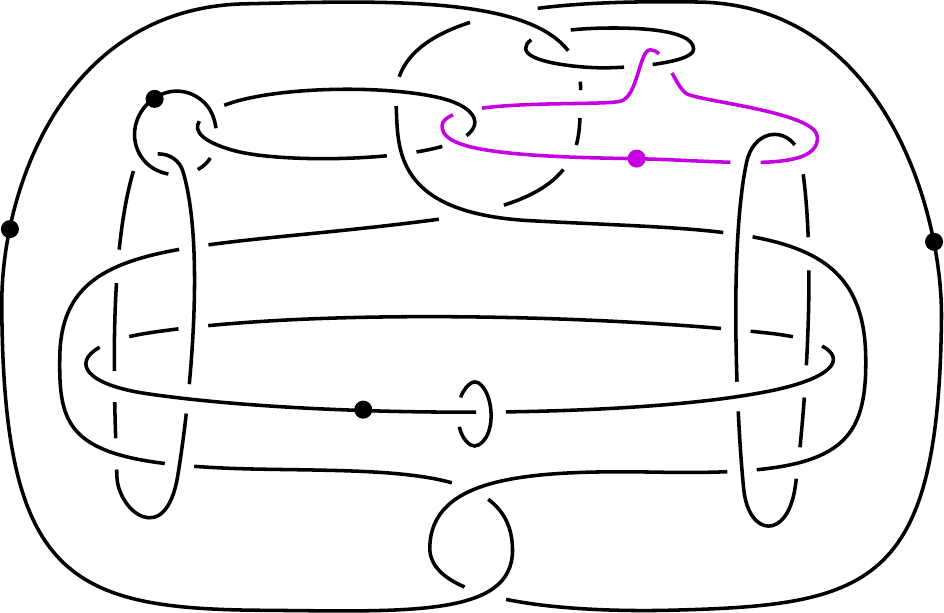}}
                            
            \put(5.62cm, 4.53cm){$-1$}
            \put(2.4cm, 4.40cm){$-1$}
            \put(3.5cm, 2.0cm){$-1$}
            \put(0.65cm, 0.75cm){$0$}
            \put(6.45cm, 0.75cm){$0$}
                        
        \end{picture}        
        }
        \caption{Gompf's infinite order cork $C(1,1;-1)$.\vspace{-8pt}}
        \label{20241024-1}
    \end{figure}

    Finally, after performing a handle slide indicated by the yellow dotted arrow,  we obtain the Kirby diagram in Figure \ref{20241024-1}. This is the same diagram as Figure 4 of \cite{Gom17a}. We shall use this diagram for the calculations in Section \ref{section-original-stein}.

    \subsection{An infinite family of infinite order corks}

    Finally, we modify the construction above to obtain an infinite family of infinite order corks. 
    
    Let $P$ be a slice pattern in the solid torus; that is, $P(U)$, the satellite of the unknot with pattern $P$, is a slice knot. Instead of drilling out the core disk $\Delta$ of the $(-1)$-framed 2-handle attached along $\{+1\}\times C_{+1}\times \{\theta_{+1}\}$ (diagrammatically, the top-right 2-handle in Figure \ref{20241024-1}), we drill out its satellite disk with pattern $P$ inside the tubular neighborhood $\nu(\Delta)$, extended to the boundary by the satellite of the annulus $I\times C_{+1}\times \{\theta_{+1}\}$. In terms of Kirby diagrams, this replaces the purple 1-handle curve in Figure \ref{20241024-1} by its satellite with pattern $P$. 
    
    Observe that the induced modification on the boundary is supported inside a tubular neighborhood of $\{-1\}\times C_{+1} \times \{\theta_{+1}\}$, so it remains disjoint from the essential region of the boundary. Moreover, this operation again merely carves out properly embedded disks, so the resulting $4$-manifold still embeds in the elliptic surface $E(n)$. Hence, by the criteria from Subsection \ref{section-understand-gompf}, this modified manifold still admits an infinite order boundary twist. Finally, if the pattern $P$ has algebraic winding number $1$, the resulting manifold is acyclic. If, in addition, $P(U)$ is the unknot and the carved-out disk is the standard slice disk, the resulting manifold is contractible, and hence is an infinite order cork. We summarize the discussion with the following auxiliary lemma.

    \begin{lemma}{\label{auxiliary}}
        Let $P$ be a slice pattern in the solid torus, and let $C_P$ be the 4-manifold obtained from the Kirby diagram of $C(1,1;-1)$ in Figure \ref{20241024-1} by replacing the purple dotted curve with its satellite with pattern $P$. Then $C_P$ admits an infinite order boundary twist. Moreover, if $P$ has algebraic winding number $1$, then $C_P$ is acyclic. If additionally $P(U)=U$, so that the satellite dotted curve is an ordinary dotted unknot, then $C_P$ is contractible.
    \end{lemma}

\begin{figure}[!htbp]
                        \centering
                        \resizebox{4.95cm}{!}{
                        \begin{picture}(5.5cm,4cm)
                            \put(0cm,0cm){\includegraphics[width=5.5cm, height=4cm]{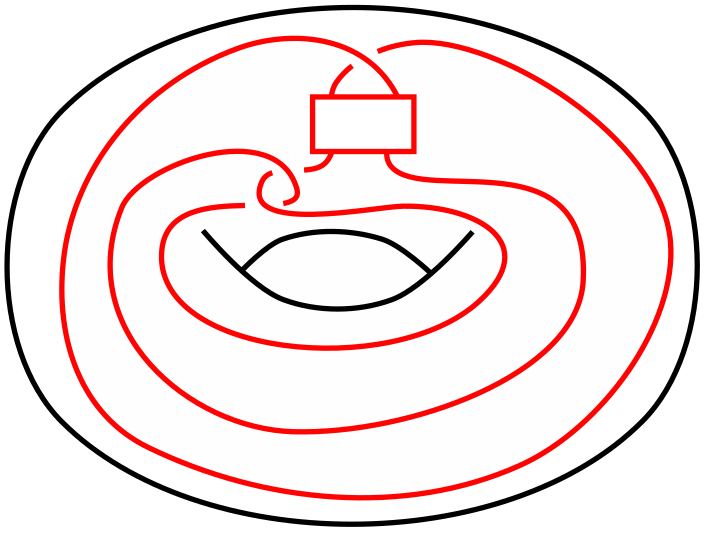}}
                            \put(2.5cm,3.0cm){\textcolor{red}{$-m$}}
                        \end{picture}
                        }
                        \caption{The satellite patterns $P_m$ used in the construction of $C_m$. \vspace{-10 pt}}
                        \label{pattern-Yasui}
                \end{figure}

    With the help of this auxiliary lemma, we define the infinite family of infinite order corks $\{C_m\}$ by replacing the purple curve with its satellite with respect to the patterns $P_m$ shown in Figure~\ref{pattern-Yasui}. To the author's knowledge, these patterns were first constructed in Yasui's paper \cite{Yas17}. The corks $C_m$ will be the starting point for the calculations in Section \ref{section-modified-stein}.

\section{Gompf's cork \texorpdfstring{$C(1,1;-1)$}{C(1,1;-1)} is Stein} {\label{section-original-stein}}

    In this section, we start with the Kirby diagram in Figure \ref{20241024-1} and prove that Gompf's cork $C(1,1;-1)$ admits a Stein structure. Notice that the black part of the diagram has the homotopy type of $S^2$. We will show that it is, in fact, diffeomorphic to a $0$-framed knot trace, i.e., the result of attaching a single $0$-framed $2$-handle to $B^4$. We will treat the purple $1$-handle as a decoration throughout the calculations until we reach the Mazur-type handle diagram on the left in Figure \ref{20240604-1}. 

\begin{figure}[H]
        \centering
        \resizebox{14.85cm}{!}{
        \begin{picture}(16.5cm,5cm)
            \put(0cm,0cm){\includegraphics[width=7.5cm, height=5cm]{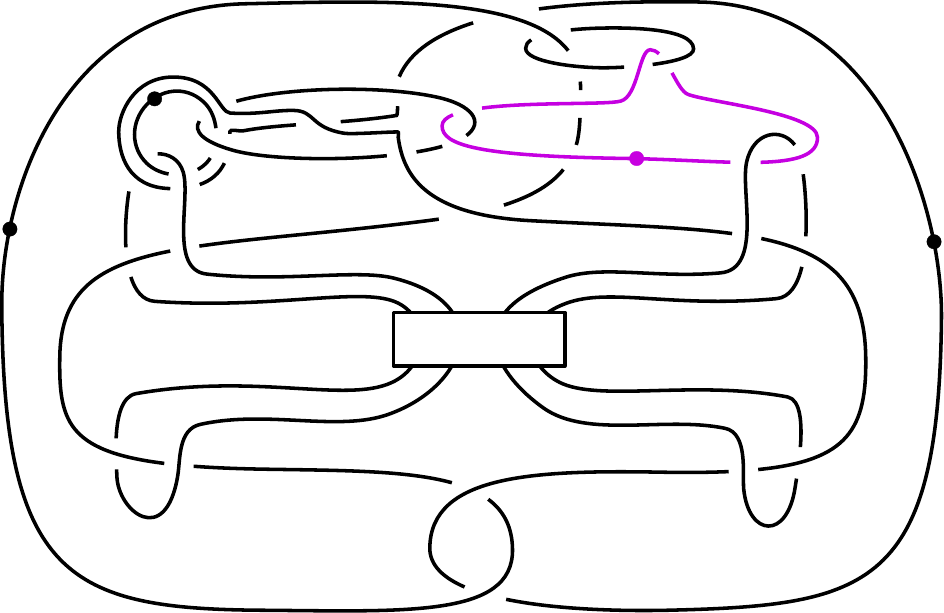}}
                            
            \put(5.62cm, 4.53cm){$-1$}
            \put(2.4cm, 4.40cm){$-1$}
            \put(3.5cm, 2.12cm){$-1$}
            \put(0.65cm, 0.75cm){$0$}
            \put(6.45cm, 0.75cm){$0$}

            \put(8.5cm,0cm){\includegraphics[width=7.5cm, height=5cm]{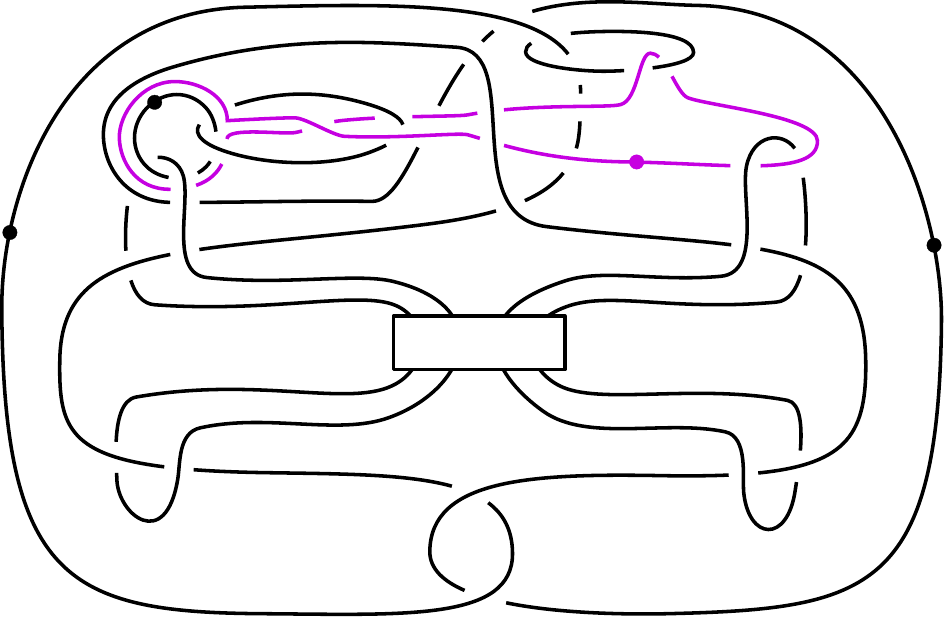}}

            \put(14.12cm, 4.53cm){$-1$}
            \put(10.9cm, 4.30cm){$-1$}
            \put(12.0cm, 2.10cm){$-1$}
            \put(9.15cm, 0.75cm){$0$}
            \put(14.95cm, 0.75cm){$0$}
                        
        \end{picture}   
        }
        \caption{\vspace{-10 pt}}
        \label{20240214-2-3}
    \end{figure}

    We start by canceling the central 1-handle with its meridional $(-1)$-framed 2-handle, followed by a 1-handle slide (with a positive half twist in the band). This gives the diagram on the left in Figure \ref{20240214-2-3}. Performing a similar handle slide for the purple curve gives us the diagram on the right.

\begin{figure}[H]
        \centering
        \resizebox{14.85cm}{!}{
        \begin{picture}(16.5cm,5cm)
            \put(0cm,0cm){\includegraphics[width=7.5cm, height=5cm]{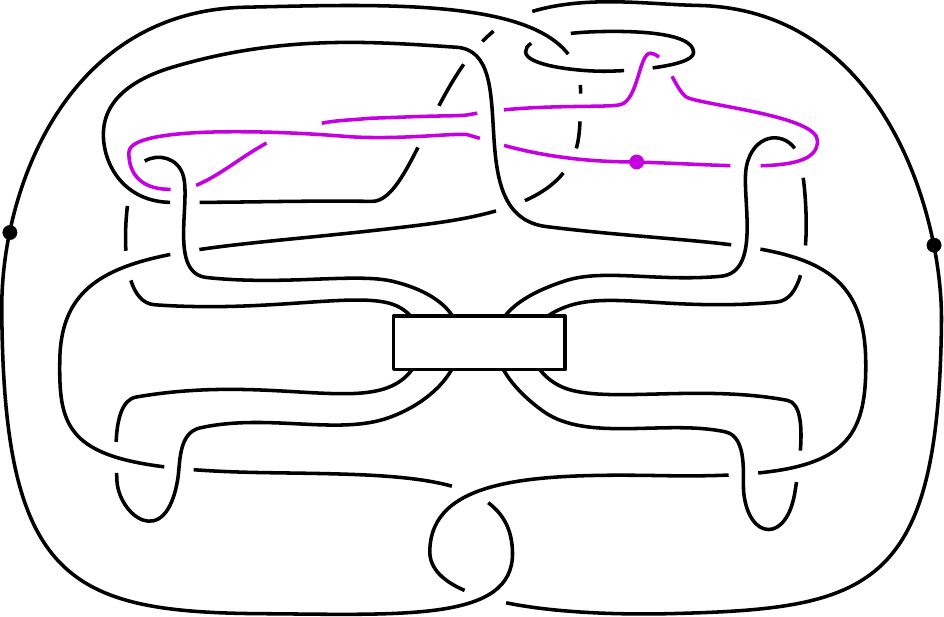}}
                            
            \put(5.62cm, 4.53cm){$-1$}
            \put(3.5cm, 2.10cm){$-1$}
            \put(0.4cm, 0.75cm){$-1$}
            \put(6.45cm, 0.75cm){$0$}

            \put(8.5cm,0cm){\includegraphics[width=7.5cm, height=5cm]{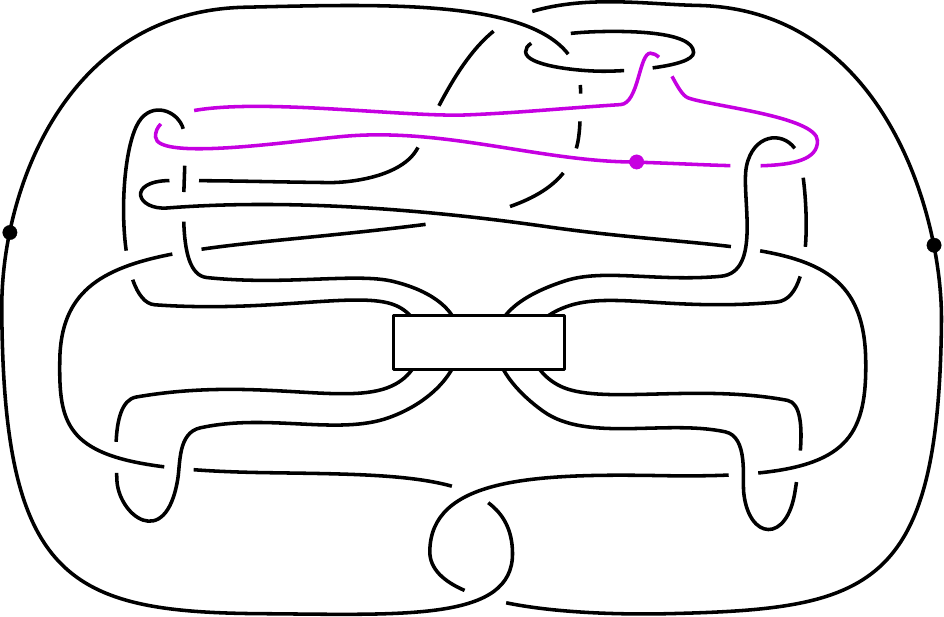}}

            \put(14.12cm, 4.53cm){$-1$}
            \put(12.0cm, 2.10cm){$-1$}
            \put(8.9cm, 0.75cm){$-1$}
            \put(14.95cm, 0.75cm){$0$}
                        
        \end{picture}        
        }
        \caption{}
        \label{20240214-4-5}
    \end{figure}

    Next, we cancel the top-left $1/2$-handle pair (after sliding the vertical $0$-framed $2$-handle over the $(-1)$-framed circle). This gives the diagram on the left in Figure \ref{20240214-4-5}. An isotopy gives the diagram on the right.

    We next aim to cancel the large 1-handle on the left with its meridional $(-1)$-framed 2-handle on the top right, namely the small one that also hooks around the purple curve. To do this, we perform a sequence of isotopies so that the two large $1$-handles appear untangled in the diagram.

\begin{figure}
        \centering
        \resizebox{14.85cm}{!}{
        \begin{picture}(16.5cm,6cm)
            \put(0cm,0cm){\includegraphics[width=7.5cm, height=6cm]{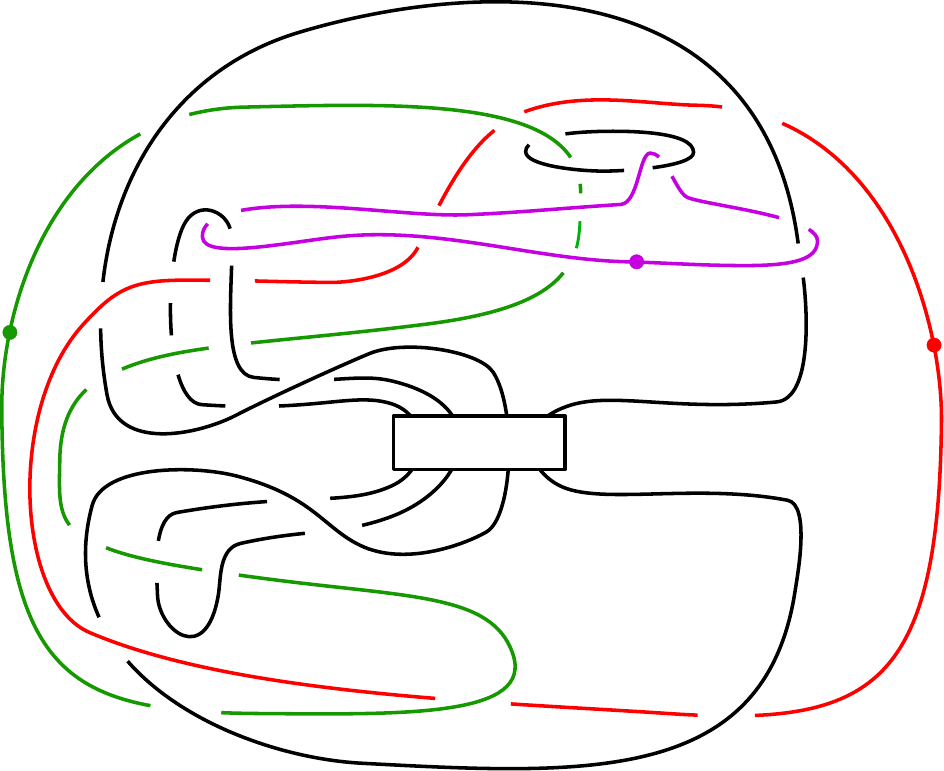}}
                            
            \put(5.53cm, 4.53cm){$-1$}
            \put(3.5cm, 2.44cm){$-1$}
            \put(1.75cm, 1.0cm){$-1$}
            \put(6.35cm, 0.75cm){$0$}

            \put(8.5cm,0cm){\includegraphics[width=7.5cm, height=6cm]{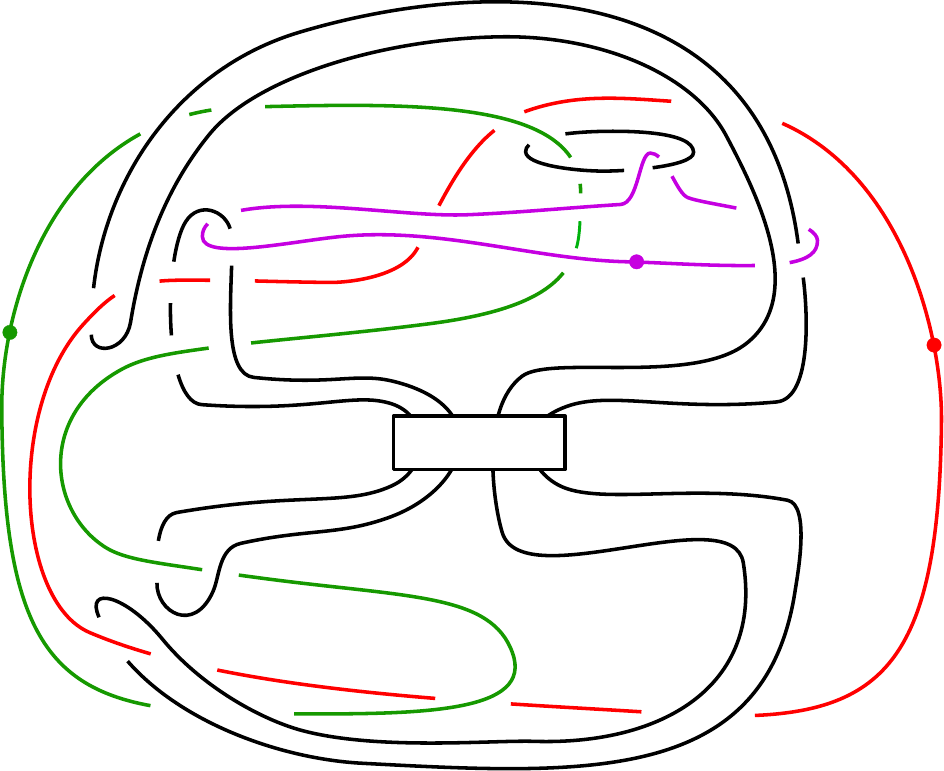}}

            \put(12.19cm, 4.53cm){$-1$}
            \put(12.0cm, 2.44cm){$-1$}
            \put(10.25cm, 1.08cm){$-1$}
            \put(14.95cm, 0.75cm){$0$}
                        
        \end{picture} 
        }
        \caption{}
        \label{20240214-6-7}
    \end{figure}

    We start by pulling the third strand of the twist box (both top and bottom) to unhook it from the large $1$-handle on the right (the red strand in Figure \ref{20240214-6-7}). This allows us to move the red strand to the left, as shown in the diagram on the left in Figure \ref{20240214-6-7}. We then pull the third strand (both top and bottom) back to obtain the diagram on the right.

\begin{figure}
        \centering
        \resizebox{13.95cm}{!}{
        \begin{picture}(15.5cm,6cm)
            \put(0cm,0cm){\includegraphics[width=7.5cm, height=6cm]{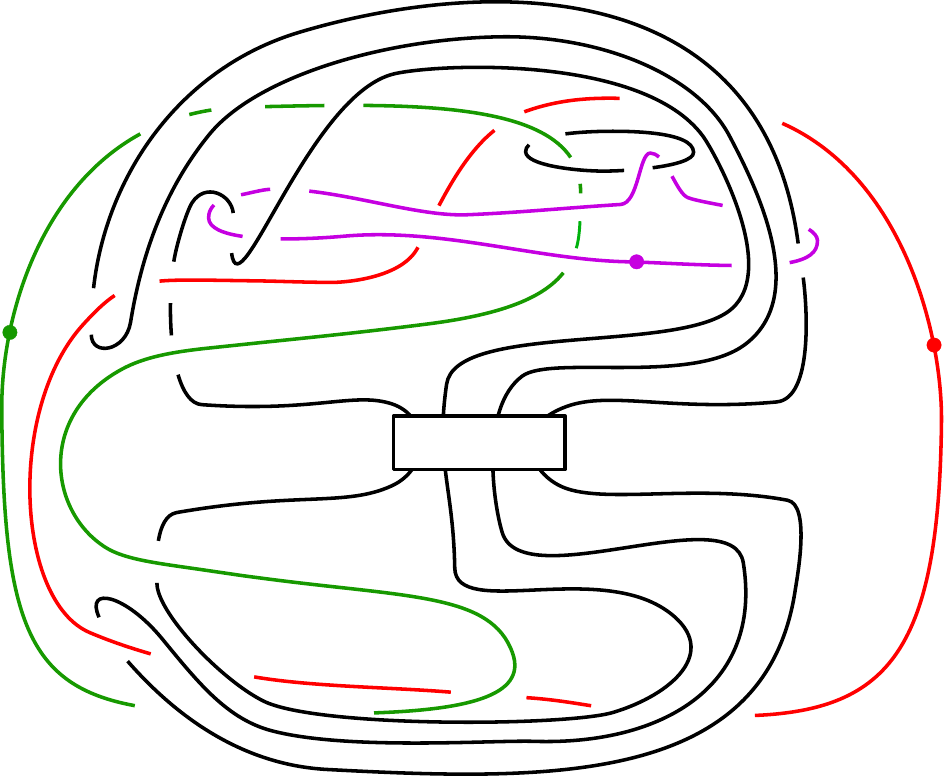}}
                            
            \put(3.72cm, 4.52cm){$-1$}
            \put(3.5cm, 2.45cm){$-1$}
            \put(1.75cm, 1.0cm){$-1$}
            \put(6.35cm, 0.75cm){$0$}

            \put(8.5cm,0cm){\includegraphics[width=6.5cm, height=6cm]{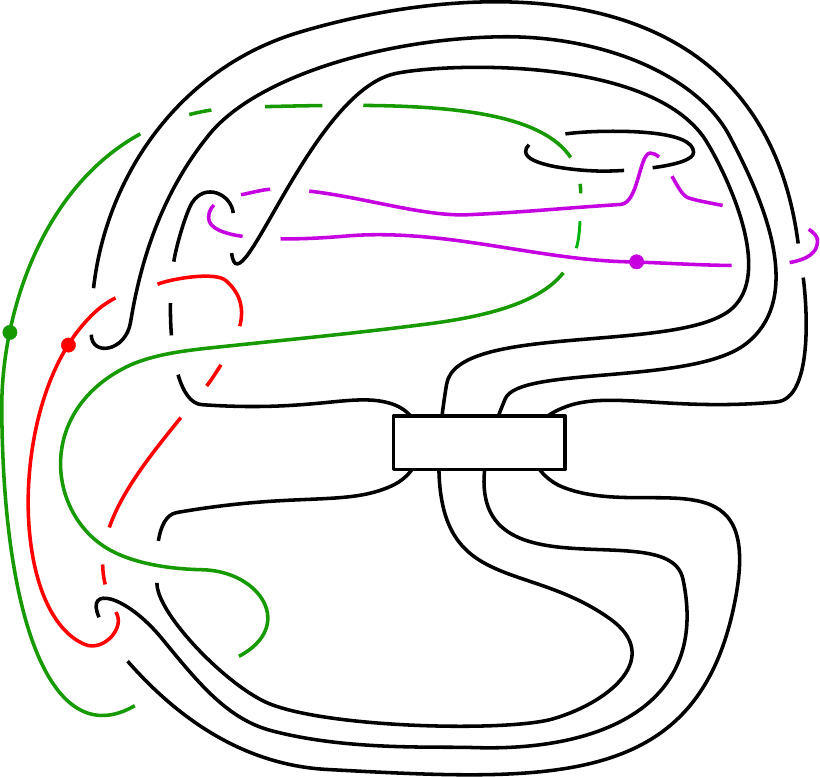}}

            \put(12.19cm, 4.53cm){$-1$}
            \put(12.05cm, 2.45cm){$-1$}
            \put(10.25cm, 1.68cm){$-1$}
            \put(14.45cm, 0.75cm){$0$}
                        
        \end{picture}   
        }
        \caption{}
        \label{20240214-8-9}
    \end{figure}

\begin{figure}
        \centering
        \resizebox{13.05cm}{!}{
        \begin{picture}(14.5cm,6cm)
            \put(0cm,0cm){\includegraphics[width=6.5cm, height=6cm]{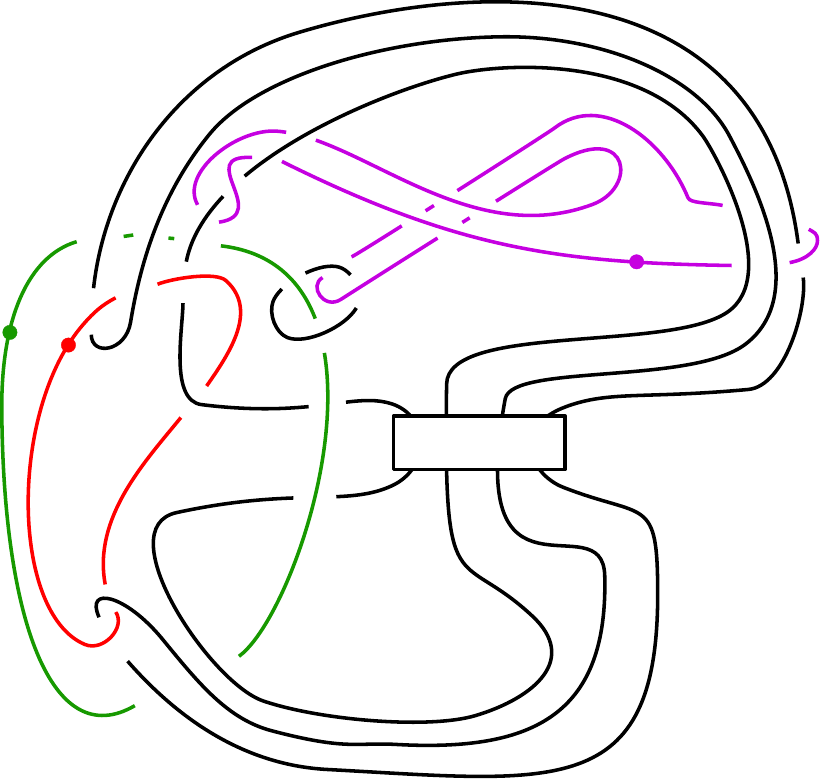}}
                            
            \put(2.92cm, 3.52cm){$-1$}
            \put(3.5cm, 2.45cm){$-1$}
            \put(2.45cm, 1.68cm){$-1$}
            \put(5.35cm, 0.75cm){$0$}

            \put(7.5cm,0cm){\includegraphics[width=6.5cm, height=6cm]{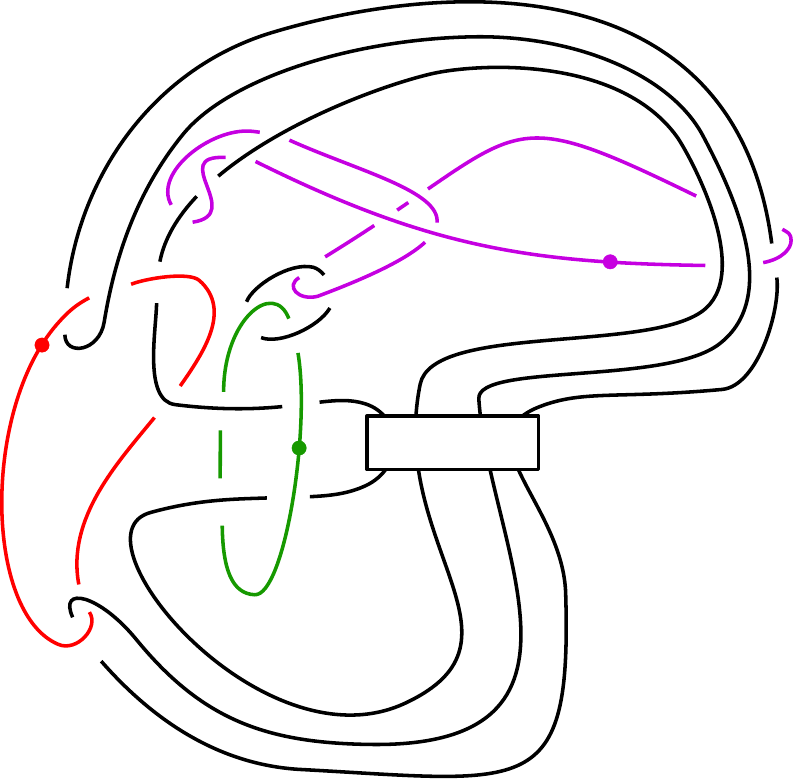}}

            \put(10.42cm, 3.52cm){$-1$}
            \put(10.95cm, 2.45cm){$-1$}
            \put(9.95cm, 1.68cm){$-1$}
            \put(12.45cm, 0.75cm){$0$}
                        
        \end{picture}
        }
        \caption{}
        \label{20240214-10-11-12}
    \end{figure}

\begin{figure}
        \centering
        \resizebox{15.75cm}{!}{
        \begin{picture}(17.5cm,6cm)
            \put(0cm,0cm){\includegraphics[width=5.5cm, height=6cm]{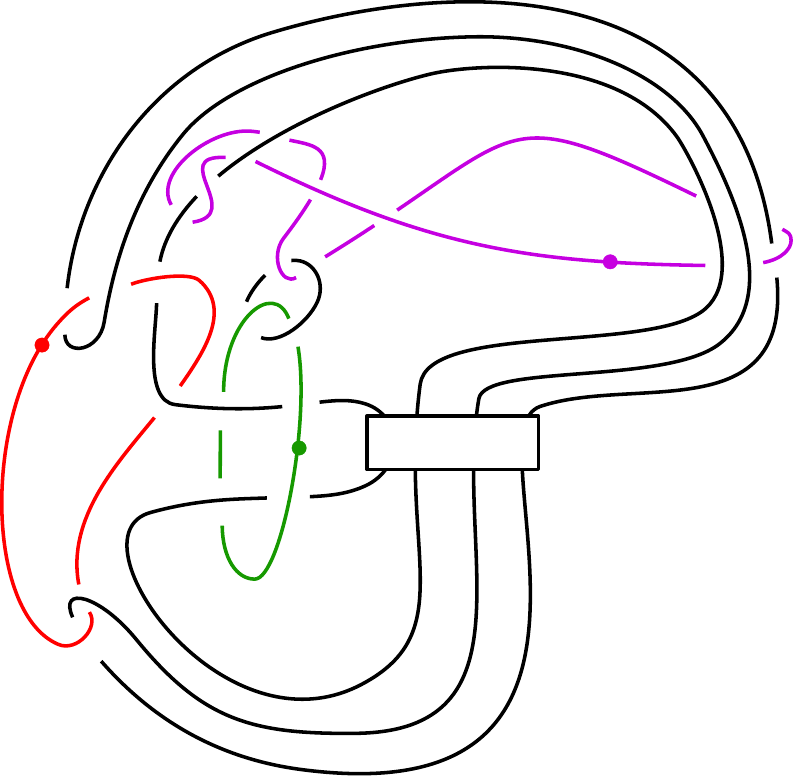}}
                            
            \put(2.30cm, 3.52cm){$-1$}
            \put(2.87cm, 2.45cm){$-1$}
            \put(2.25cm, 1.78cm){$-1$}
            \put(3.8cm, 0.75cm){$0$}

            \put(6.0cm,0cm){\includegraphics[width=5.5cm, height=6cm]{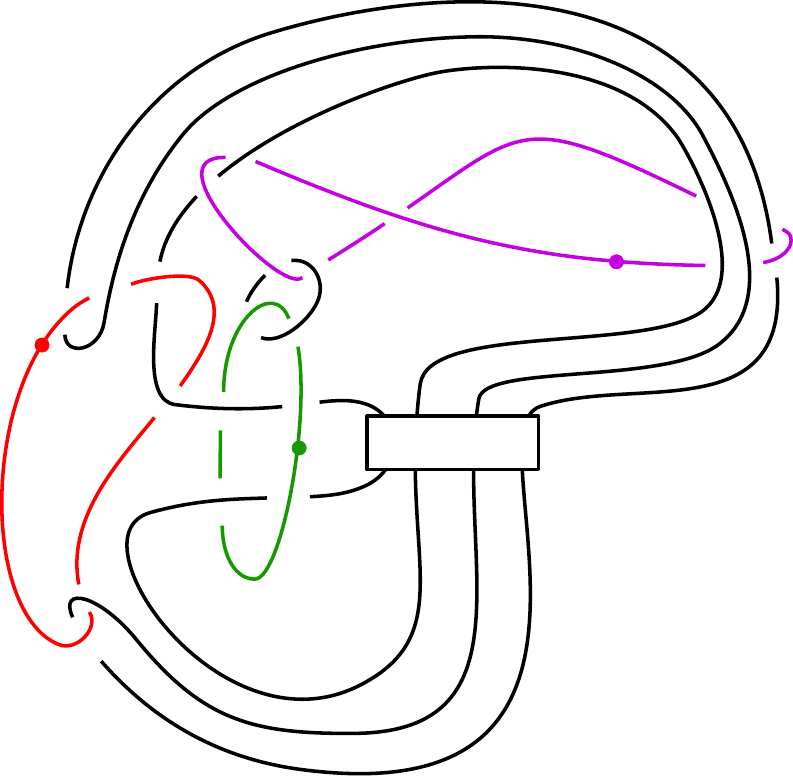}}

            \put(8.30cm, 3.52cm){$-1$}
            \put(8.87cm, 2.45cm){$-1$}
            \put(8.25cm, 1.78cm){$-1$}
            \put(9.8cm, 0.75cm){$0$}

            \put(12.0cm,0cm){\includegraphics[width=5.5cm, height=6cm]{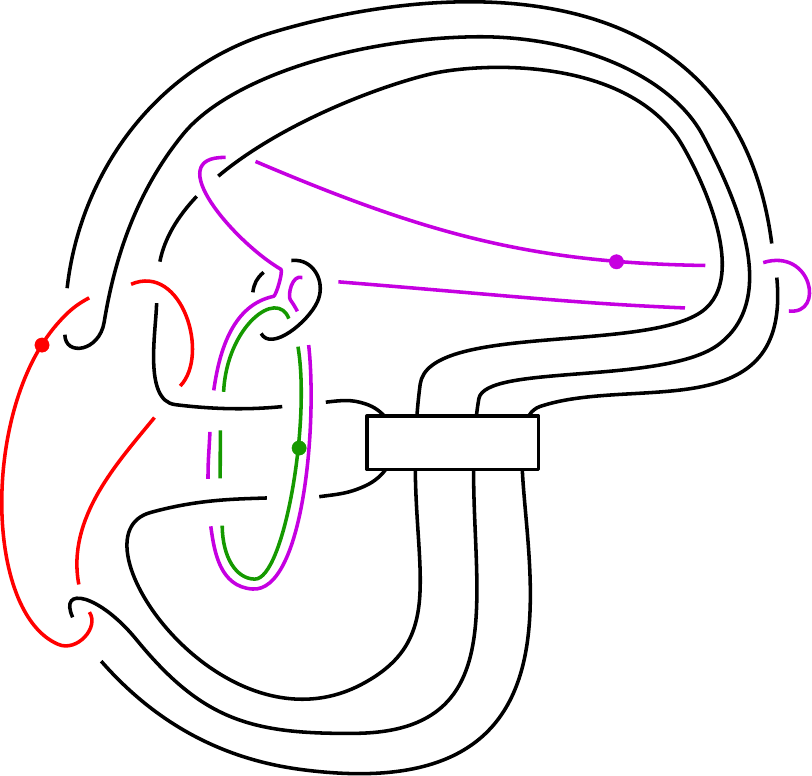}}
            \put(14.30cm, 3.52cm){$-1$}
            \put(14.87cm, 2.45cm){$-1$}
            \put(14.25cm, 1.78cm){$-1$}
            \put(15.8cm, 0.75cm){$0$}
                        
        \end{picture}
        }
        \caption{}
        \label{20240214-13-14-15}
    \end{figure}

    Next, we pull the second strand (again both top and bottom) of the twist box. This gives us the picture on the left in Figure \ref{20240214-8-9}. We then carefully pull the red and green 1-handles away from each other. This process, shown in Figures \ref{20240214-8-9}, \ref{20240214-10-11-12} and \ref{20240214-13-14-15}, eventually gives us the diagram in the middle of Figure \ref{20240214-13-14-15}. If we slide the purple curve over the green 1-handle, we can put the green 1-handle and its meridional $2$-handle in canceling position. This gives the last diagram in Figure \ref{20240214-13-14-15}.

\begin{figure}
        \centering
        \resizebox{15.75cm}{!}{
        \begin{picture}(17.5cm,6cm)
            \put(0cm,0cm){\includegraphics[width=5.5cm, height=6cm]{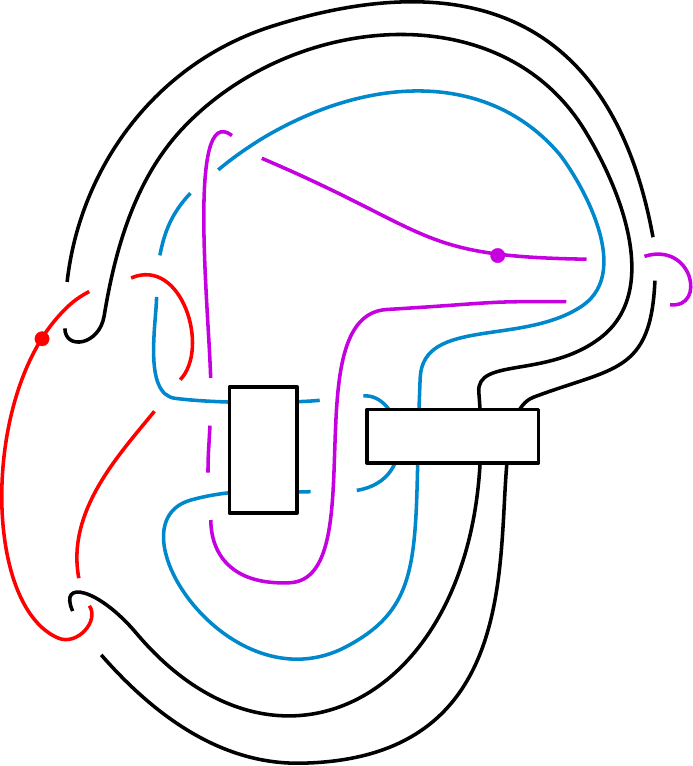}}
                            
            \put(1.85cm, 2.4cm){$-1$}
            \put(3.31cm, 2.45cm){$-1$}
            \put(2.75cm, 1.78cm){\textcolor{blue}{$-1$}}
            \put(4.0cm, 0.75cm){$0$}

            \put(6.0cm,0cm){\includegraphics[width=5.5cm, height=6cm]{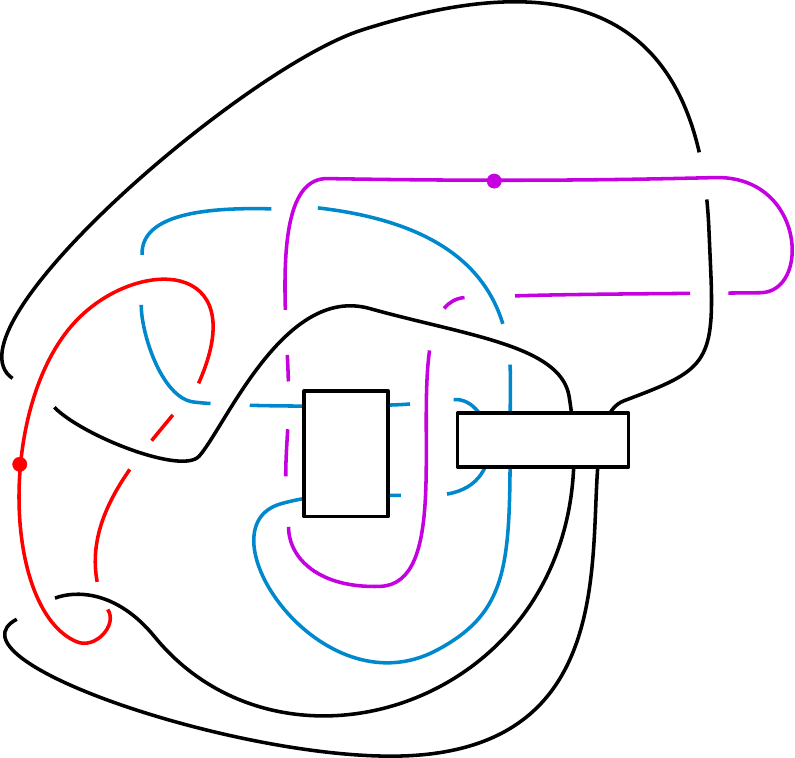}}

            \put(8.15cm, 2.3cm){$-1$}
            \put(9.5cm, 2.4cm){$-1$}
            \put(8.98cm, 1.7cm){\textcolor{blue}{$-1$}}
            \put(10.22cm, 0.75cm){$0$}

            \put(12.0cm,0cm){\includegraphics[width=5.5cm, height=6cm]{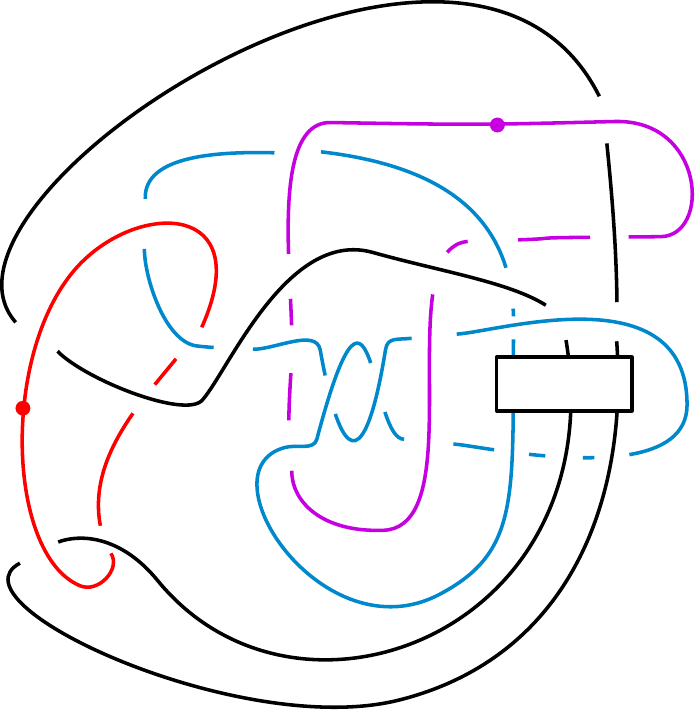}}
            
            \put(16.21cm, 2.63cm){$-1$}
            \put(14.55cm, 1.89cm){\textcolor{blue}{$-1$}}
            \put(16.6cm, 0.75cm){$0$}
                        
        \end{picture}
        }
        \caption{}
        \label{20240214-16-17-18}
    \end{figure}

\begin{figure}[H]
        \centering
        \resizebox{15.75cm}{!}{
        \begin{picture}(17.5cm,6cm)
            \put(0cm,0cm){\includegraphics[width=5.5cm, height=6cm]{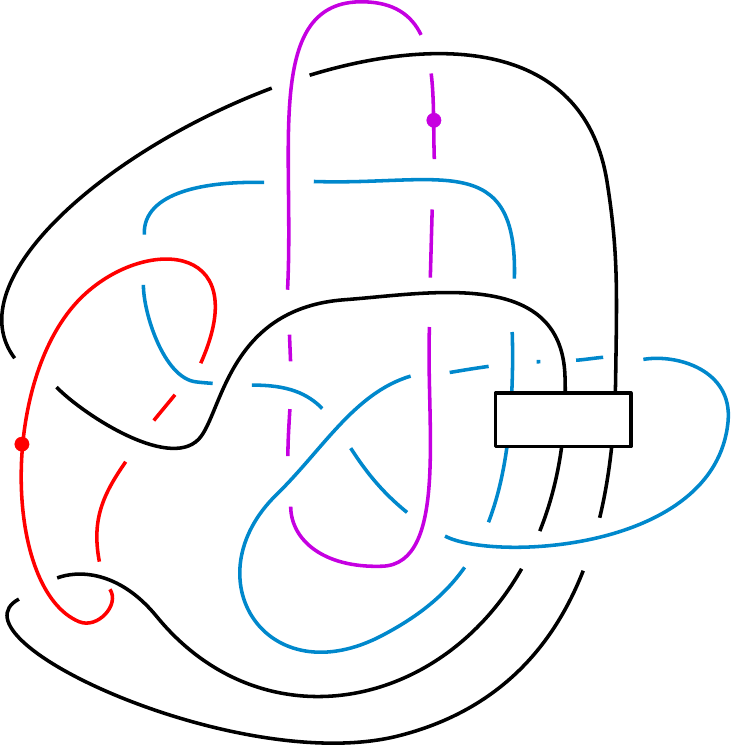}}
                            
            \put(1.55cm, 2.35cm){$-1$}
            \put(3.95cm, 2.47cm){$-1$}
            \put(2.2cm, 0.95cm){\textcolor{blue}{$-1$}}
            \put(4.3cm, 0.75cm){$0$}

            \put(6.0cm,0cm){\includegraphics[width=5.5cm, height=6cm]{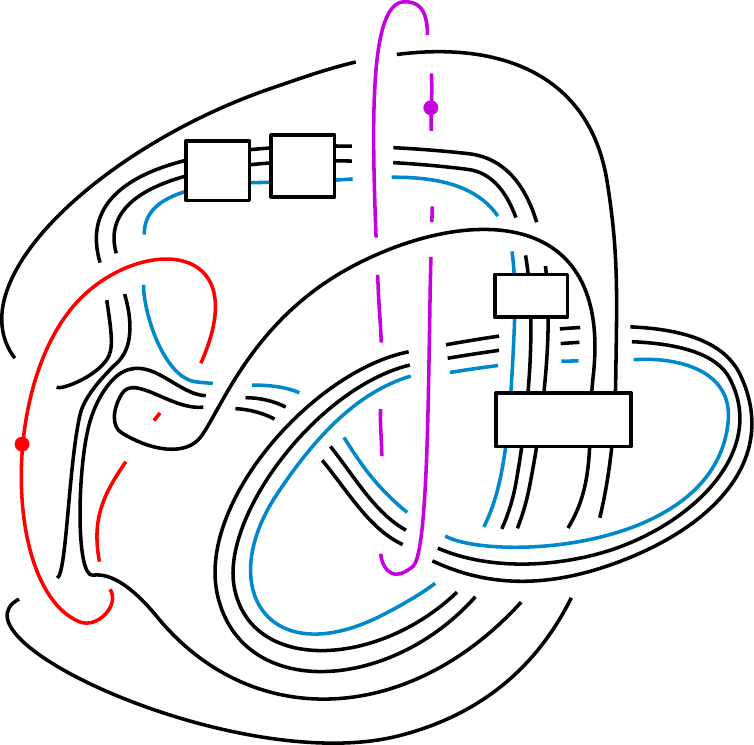}}

            \put(8.0cm, 4.56cm){\footnotesize{$-1$}}
            \put(7.37cm, 4.53cm){\footnotesize{$-3$}}
            \put(9.86cm, 2.5cm){\small{$-1$}}
            \put(9.66cm, 3.52cm){\small{\footnotesize{$+1$}}}
            \put(7.98cm, 1.1cm){\textcolor{blue}{$-1$}}
            \put(10.22cm, 0.75cm){$0$}

            \put(12.0cm,0cm){\includegraphics[width=5.5cm, height=6cm]{plots-final/Figure-26-3.pdf}}
            \put(16.22cm, 0.75cm){$0$}
            \put(15.87cm, 3.05cm){$-1$}
            \put(15.17cm, 5.43cm){\footnotesize{$-3$}}
                        
        \end{picture}
        }
        \caption{}
        \label{20240214-19-20-21}
    \end{figure}

    Now we cancel the green $1$-handle with its meridional $2$-handle. This gives the diagram on the left in Figure \ref{20240214-16-17-18}. Notice that the attaching circle of the $(-1)$-framed 2-handle had linking number $0$ with the green $1$-handle, so its framing remains unchanged. An isotopy gives the rightmost diagram. A further isotopy gives the first picture of Figure \ref{20240214-19-20-21}.

    To prepare for the eventual handle cancellation, we (double) slide the black $2$-handle over the blue $2$-handle. This gives the middle diagram in Figure \ref{20240214-19-20-21}. Note that the $(-3)$-twist box compensates for the writhe of the blue $2$-handle (i.e., right-handed trefoil), while the $(-1)$-twist box comes from the $(-1)$-framing of the target curve. Also note that the $(+1)$-twist box compensates for the handle slide passing through the large $(-1)$-twist box. Finally, since the two slides use oppositely oriented parallels of the blue curve, their framing contributions cancel. Hence, the framing of the slid handle remains unchanged.

    Finally, we cancel the blue/red 1/2-handle pair by directly erasing them. What is left is the rightmost diagram in Figure \ref{20240214-19-20-21}, which is exactly the left-hand diagram in Figure \ref{20240604-1}.

\begin{figure}[H]
        \centering
        \resizebox{4.59cm}{!}{
        \begin{picture}(5.1cm,5.1cm)
            \put(0cm,0cm){\includegraphics[width=5.1cm, height=5.1cm]{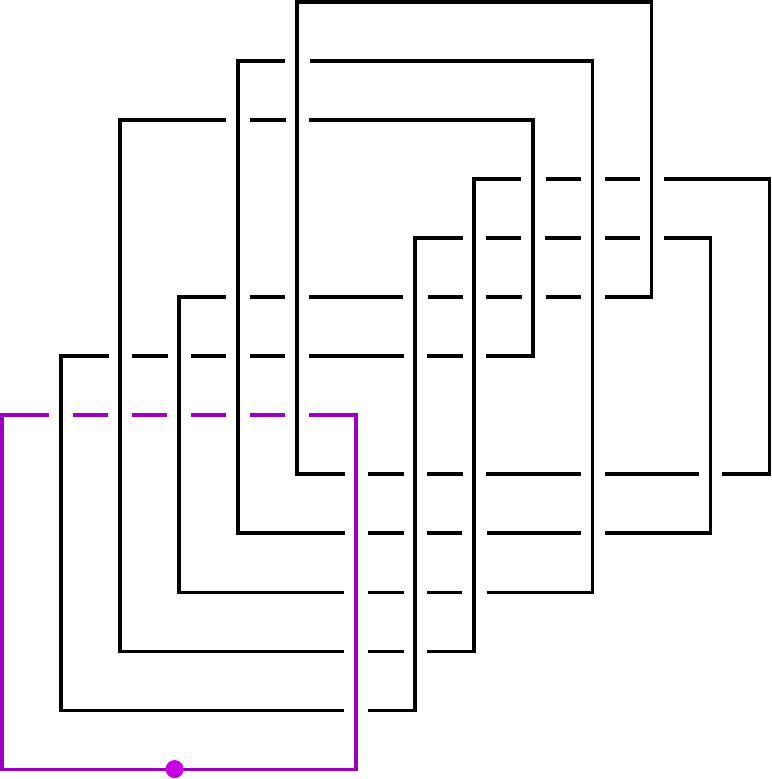}}
            \put(4.3cm, 1.2cm){$0$}
        \end{picture}
        }
        \caption{A grid diagram for the cork $C(1,1;-1).$}
        \label{20240420-1}
    \end{figure}

    Next, we wish to write down a Legendrian representation of the cork $C(1,1;-1)$. With the help of Gridlink~\cite{GridLink}, we obtain Figure \ref{20240420-1}.

\begin{figure}[H]
        \centering
        \resizebox{15.75cm}{!}{
        \begin{picture}(17.5cm,6cm)
            \put(0cm,0cm){\includegraphics[width=5.5cm, height=6cm]{plots-final/Figure-26-3.pdf}}
                            
            \put(4.22cm, 0.75cm){$0$}
            \put(3.87cm, 3.05cm){$-1$}
            \put(3.17cm, 5.43cm){\footnotesize{$-3$}}

            \put(6.0cm,0cm){\includegraphics[width=5.5cm, height=6cm]{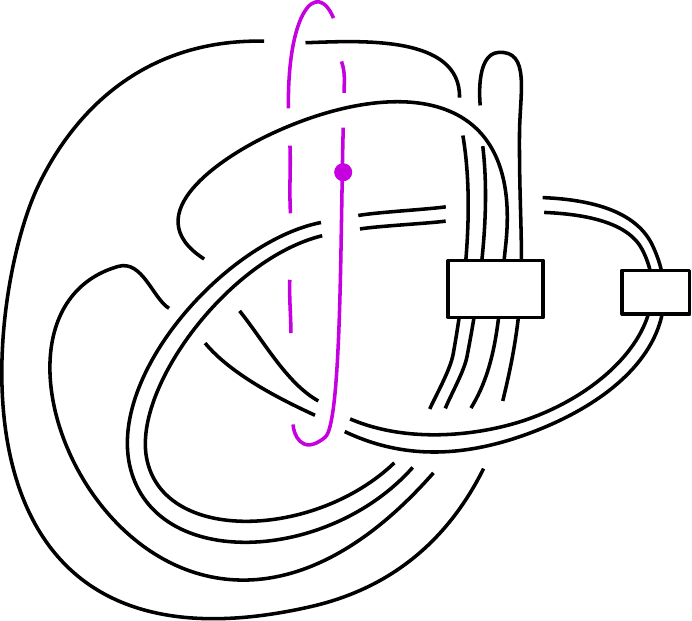}}

            \put(10.02cm, 0.75cm){$0$}
            \put(9.68cm, 3.09cm){$-1$}
            \put(10.99cm, 3.105cm){\scriptsize{$-\frac{5}{2}$}}

            \put(12.0cm,0cm){\includegraphics[width=5.5cm, height=6cm]{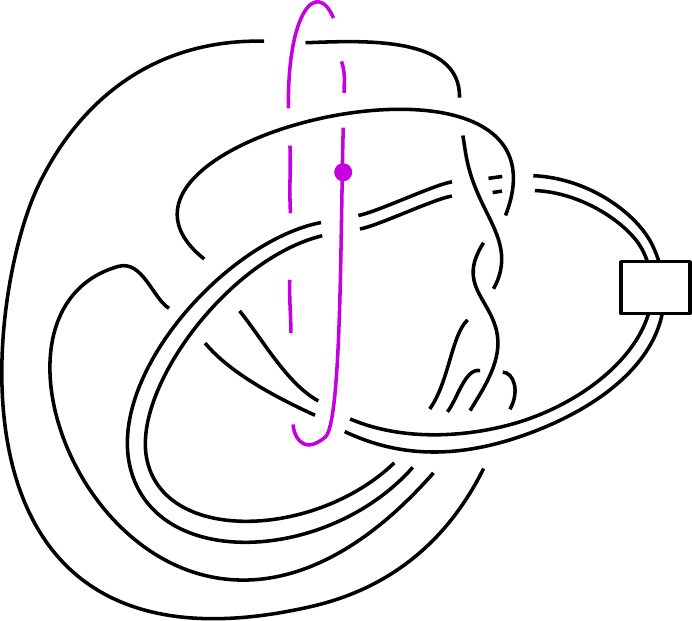}}
            \put(16.02cm, 0.75cm){$0$}
            \put(16.99cm, 3.16cm){\scriptsize{$-\frac{5}{2}$}}
                        
        \end{picture}
        }
        \caption{}
        \label{20240420-2-3-4}
    \end{figure}

\begin{figure}[H]
        \centering
        \resizebox{15.75cm}{!}{
        \begin{picture}(17.5cm,5.5cm)
            \put(0cm,0cm){\includegraphics[width=5.5cm, height=5.5cm]{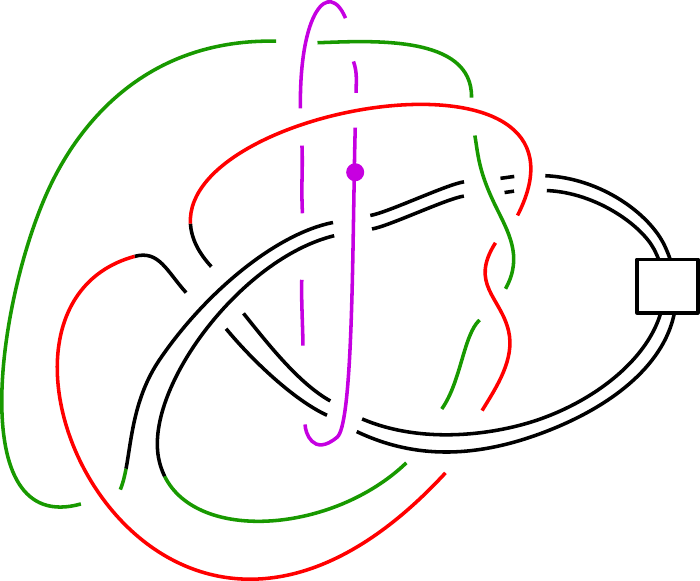}}
                            
            \put(4.22cm, 0.99cm){$0$}
            \put(5.03cm, 2.73cm){\scriptsize{$-\frac{5}{2}$}}

            \put(6.0cm,0cm){\includegraphics[width=5.5cm, height=5.5cm]{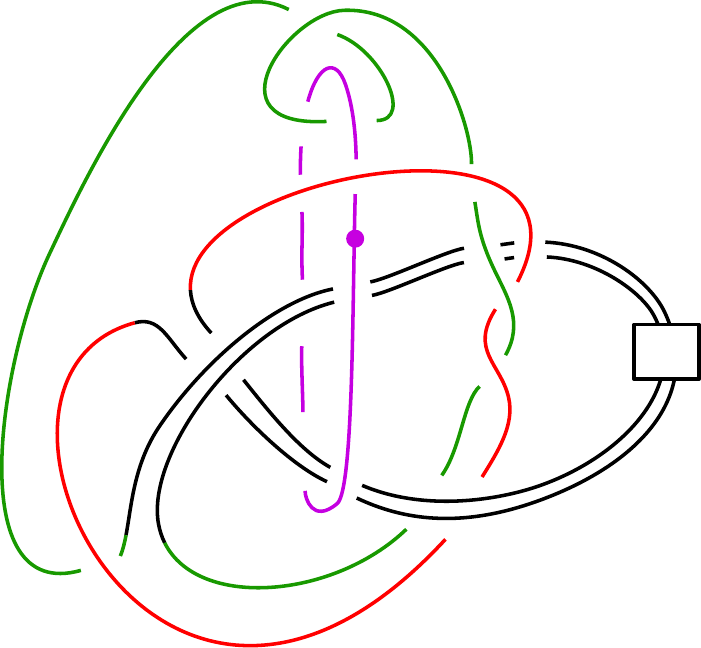}}

            \put(10.32cm, 0.95cm){$0$}
            \put(11.02cm, 2.44cm){\scriptsize{$-\frac{5}{2}$}}

            \put(12.0cm,0cm){\includegraphics[width=5.5cm, height=5.5cm]{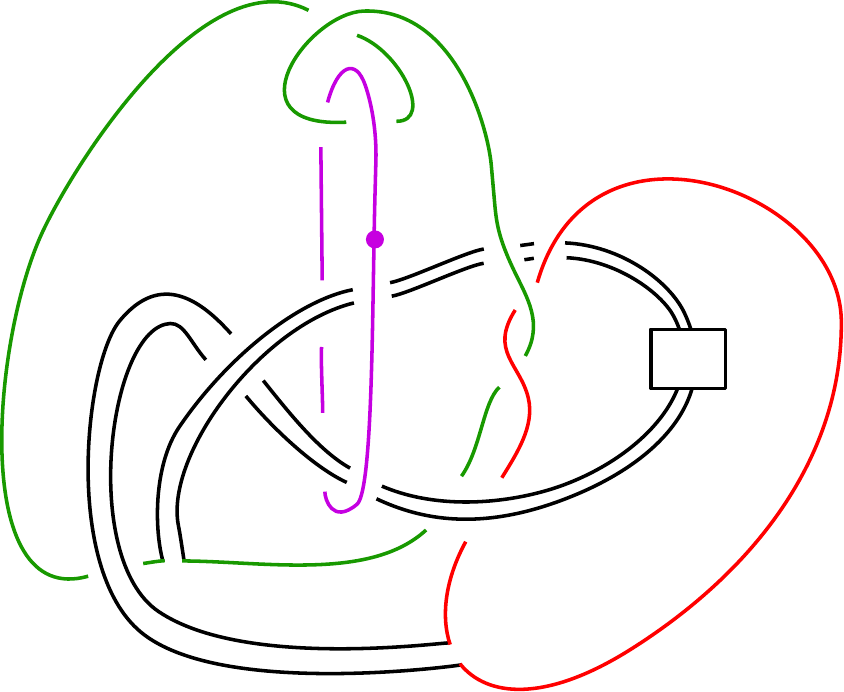}}
            \put(16.02cm, 1.35cm){$0$}
            \put(16.27cm, 2.58cm){\scriptsize{$-\frac{5}{2}$}}
                        
        \end{picture}
        }
        \caption{}
        \label{20240420-5-6-7}
    \end{figure}

\begin{figure}[H]
        \centering
        \resizebox{15.75cm}{!}{
        \begin{picture}(17.5cm,5.5cm)
            \put(0cm,0cm){\includegraphics[width=5.5cm, height=5.5cm]{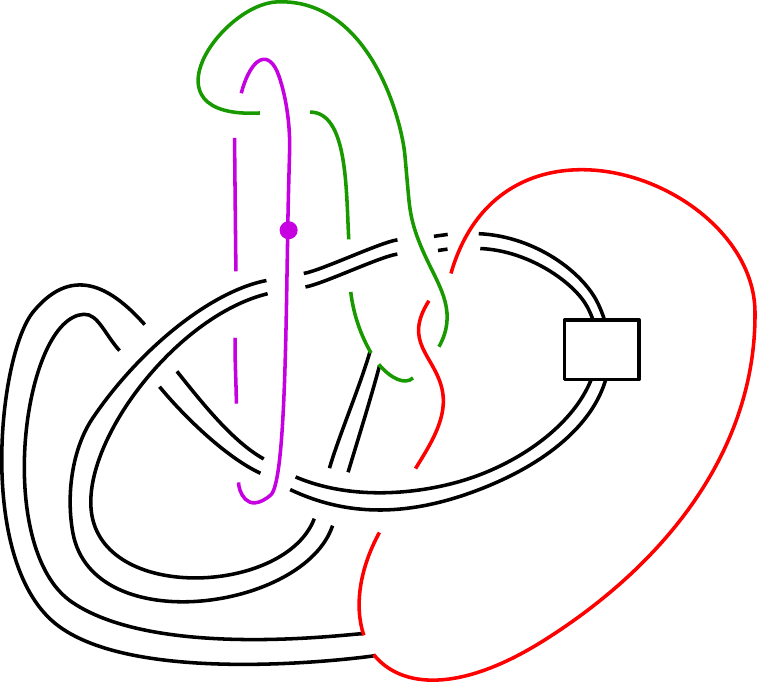}}
                            
            \put(4.92cm, 0.79cm){$0$}
            \put(4.14cm, 2.57cm){\small{$-2$}}

            \put(6.0cm,0cm){\includegraphics[width=5.5cm, height=5.5cm]{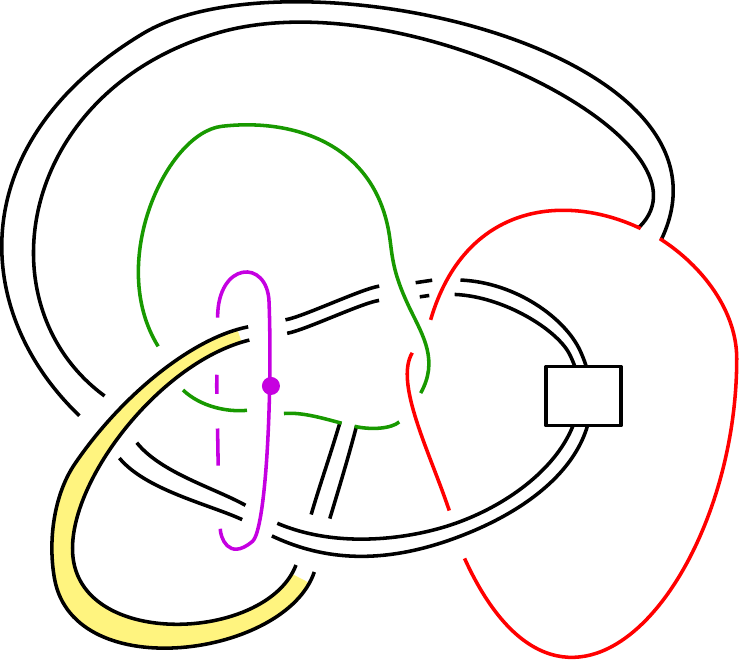}}

            \put(11.32cm, 0.75cm){$0$}
            \put(10.11cm, 2.09cm){\small{$-2$}}

            \put(12.0cm,0cm){\includegraphics[width=5.5cm, height=5.5cm]{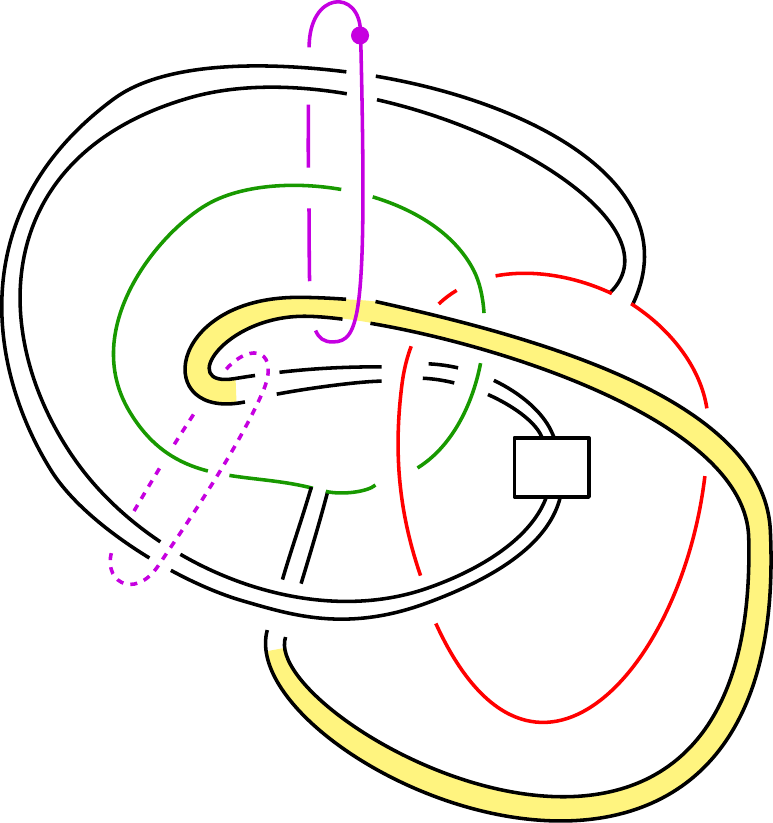}}
            \put(16.62cm, 0.75cm){$0$}
            \put(15.71cm, 2.28cm){\footnotesize{$-2$}}
                        
        \end{picture}
        }
        \caption{For the second step, first move the highlighted strip in the second diagram to its position in the last diagram, then rotate the purple circle accordingly.}
        \label{20240420-8-9-10}
    \end{figure}
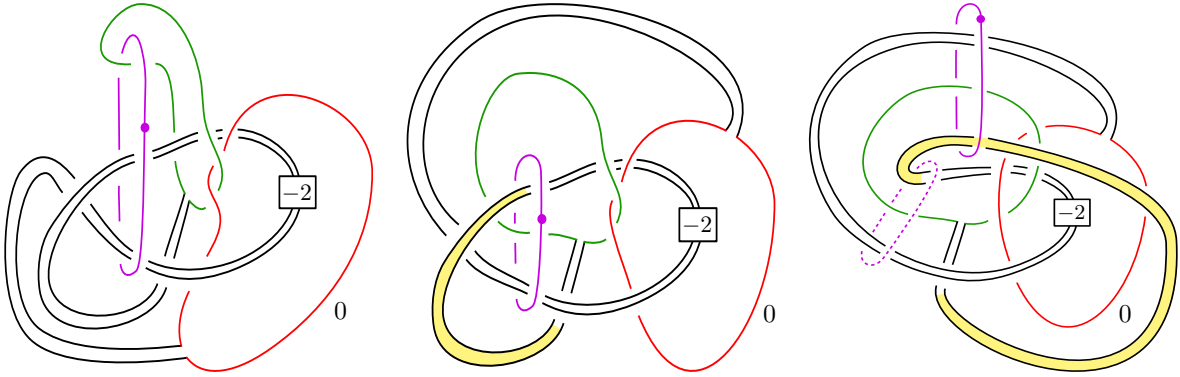

\begin{figure}[H]
        \centering
        \resizebox{15.75cm}{!}{
        \begin{picture}(17.5cm,5.5cm)
            \put(0cm,0cm){\includegraphics[width=5.5cm, height=5.5cm]{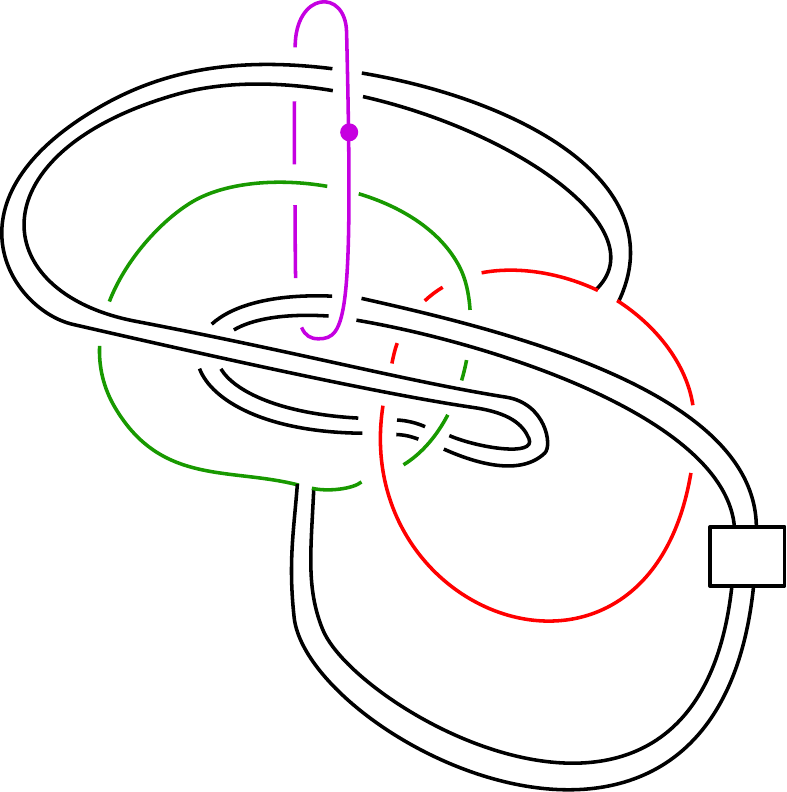}}
                            
            \put(5.35cm, 0.55cm){$0$}
            \put(4.99cm, 1.53cm){\small{$-3$}}

            \put(6.0cm,0cm){\includegraphics[width=5.5cm, height=5.5cm]{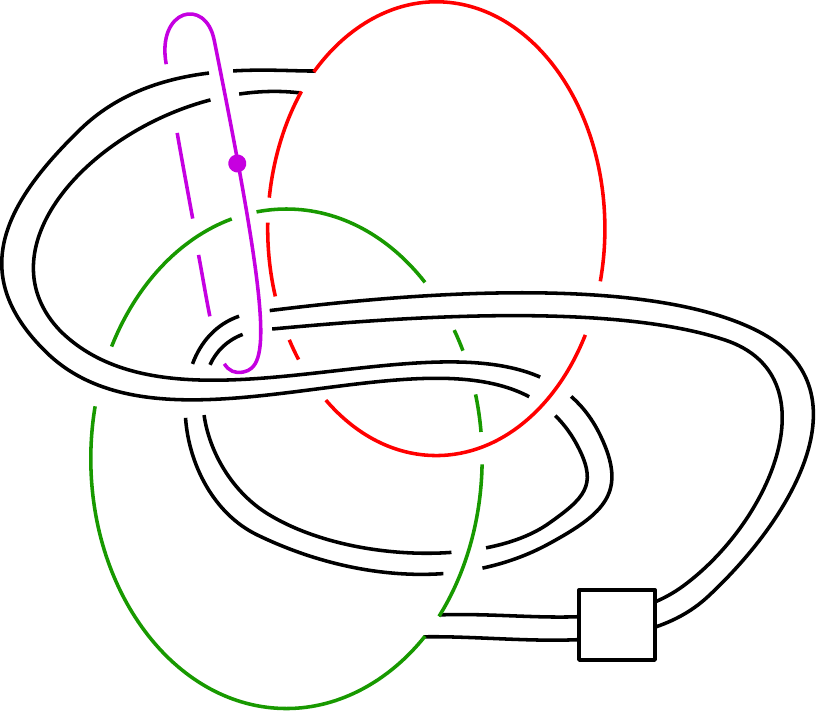}}

            \put(11.12cm, 0.75cm){$0$}
            \put(9.92cm, 0.54cm){\small{$-3$}}

            \put(12.0cm,0cm){\includegraphics[width=5.5cm, height=5.5cm]{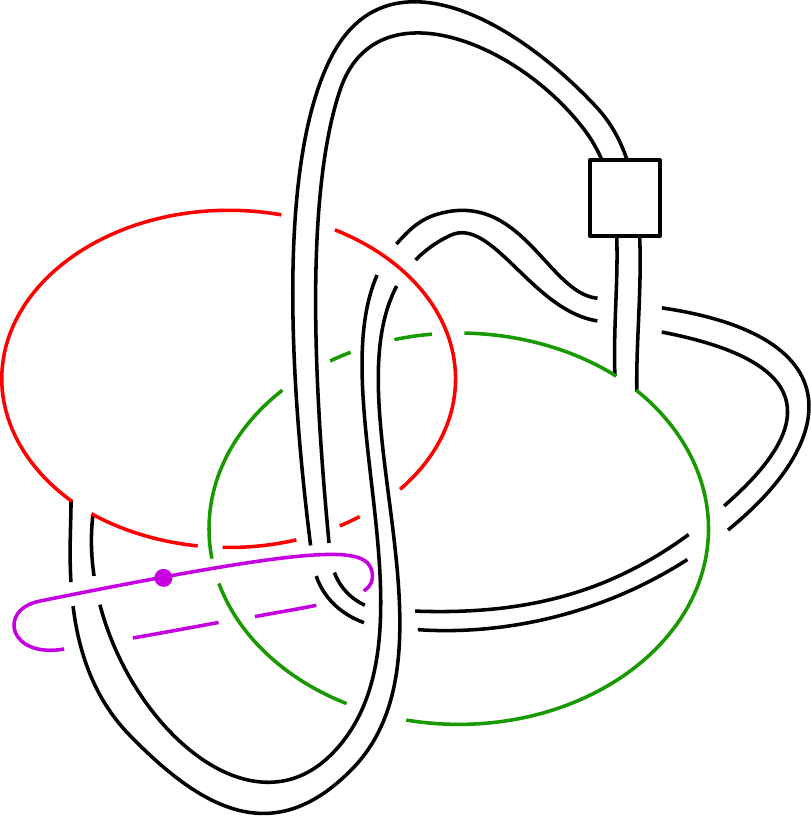}}
            \put(17.12cm, 1.65cm){$0$}
            \put(16.02cm, 4.06cm){\footnotesize{$-3$}}
                        
        \end{picture}
        }
        \caption{}
        \label{20240420-11-12-13}
    \end{figure}

\begin{figure}[H]
        \centering
        \resizebox{15.3cm}{!}{
        \begin{picture}(17cm,5cm)
            \put(0cm,0cm){\includegraphics[width=5cm, height=5cm]{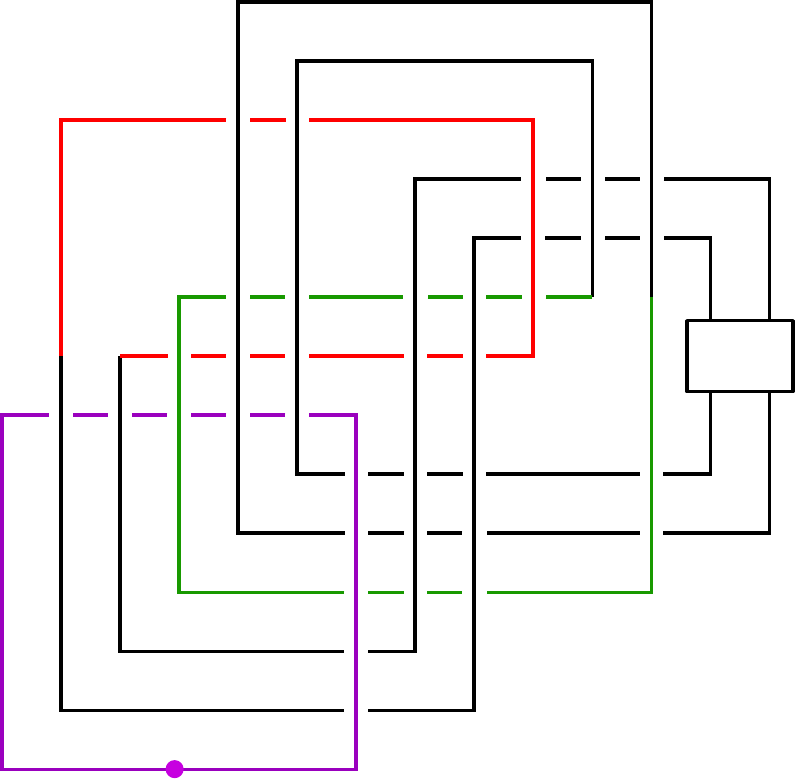}}
                            
            \put(3.45cm, 0.75cm){$0$}
            \put(4.43cm, 2.62cm){\small{$-3$}}

            \put(6.0cm,0cm){\includegraphics[width=5cm, height=5cm]{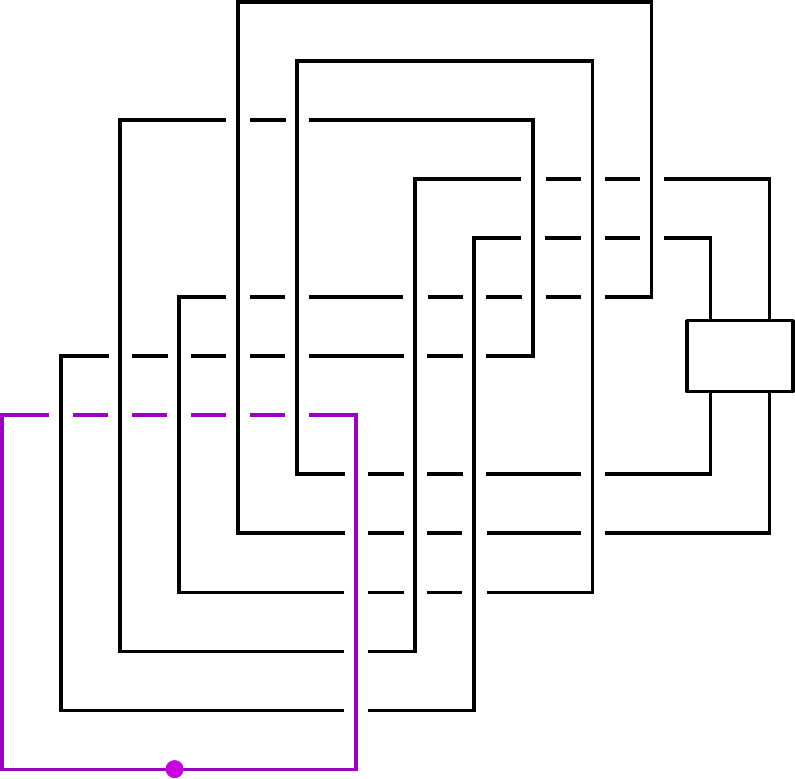}}

            \put(10.3cm, 1.12cm){$0$}
            \put(10.43cm, 2.62cm){\small{$-2$}}

            \put(12.0cm,0cm){\includegraphics[width=5cm, height=5cm]{plots-final/Figure-27.pdf}}
            \put(16.2cm, 1.12cm){$0$}
                        
        \end{picture}
        }
        \caption{\vspace{-10 pt}}
        \label{20240420-16-17-18}
    \end{figure}

    The fact that the link in Figure \ref{20240420-1} is isotopic to the last link in Figure \ref{20240214-19-20-21} is shown in the sequence of diagrams from Figure \ref{20240420-2-3-4} through Figure \ref{20240420-16-17-18}. 

\begin{figure}[H]
                        \centering
                        \resizebox{13.5cm}{!}{
                        \begin{picture}(15cm,5cm)
                            \put(0cm,0.0cm){\includegraphics[width=7cm, height=5cm]{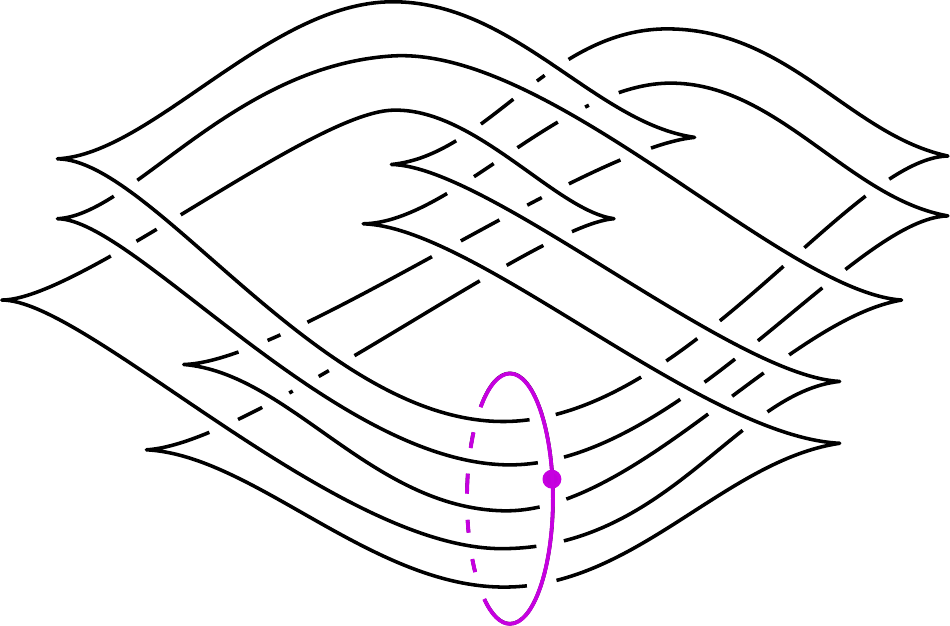}}
                            \put(8cm,0cm){\includegraphics[width=7cm, height=5cm]{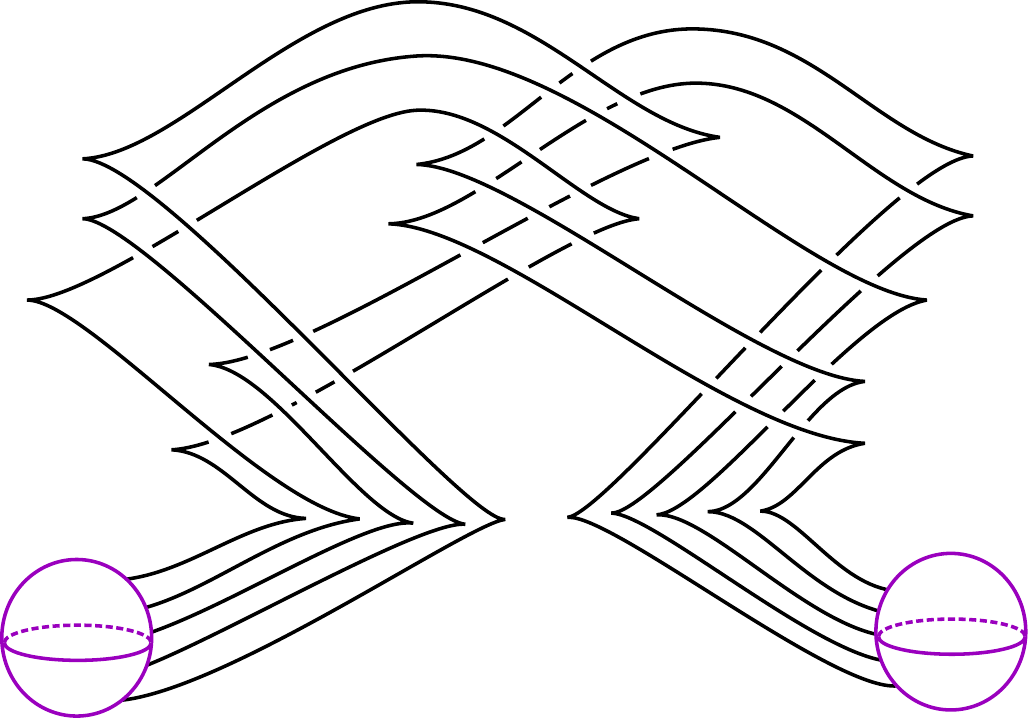}}

                            \put(6.35cm,4.32cm){$0$}
                            \put(14.35cm,4.32cm){$0$}
                        
                        \end{picture}   
                        }
                        \caption{A Legendrian diagram for the cork $C(1,1;-1)$.\vspace{-10 pt}}
                        \label{20240420-19-20}
            \end{figure}

    On the other hand, a grid diagram determines a Legendrian diagram, obtained by a rotation of $45^\circ$ counterclockwise. Thus, we obtain the first diagram in Figure \ref{20240420-19-20}. At this stage, the Thurston-Bennequin number of the black curve is $+1$.

    Next, we put the Legendrian diagram into standard position; that is, we regard the black curve as a Legendrian knot in $S^1\times S^2$. This yields the second diagram in Figure \ref{20240420-19-20}. At this stage, the Thurston-Bennequin number of the black curve is $-4$.

\begin{figure}[H]
                        \centering
                        \resizebox{4.77cm}{!}{
                        \begin{picture}(5.3cm,5.5cm)
                            \put(0cm,0cm){\includegraphics[width=6.3cm, height=5.5cm]{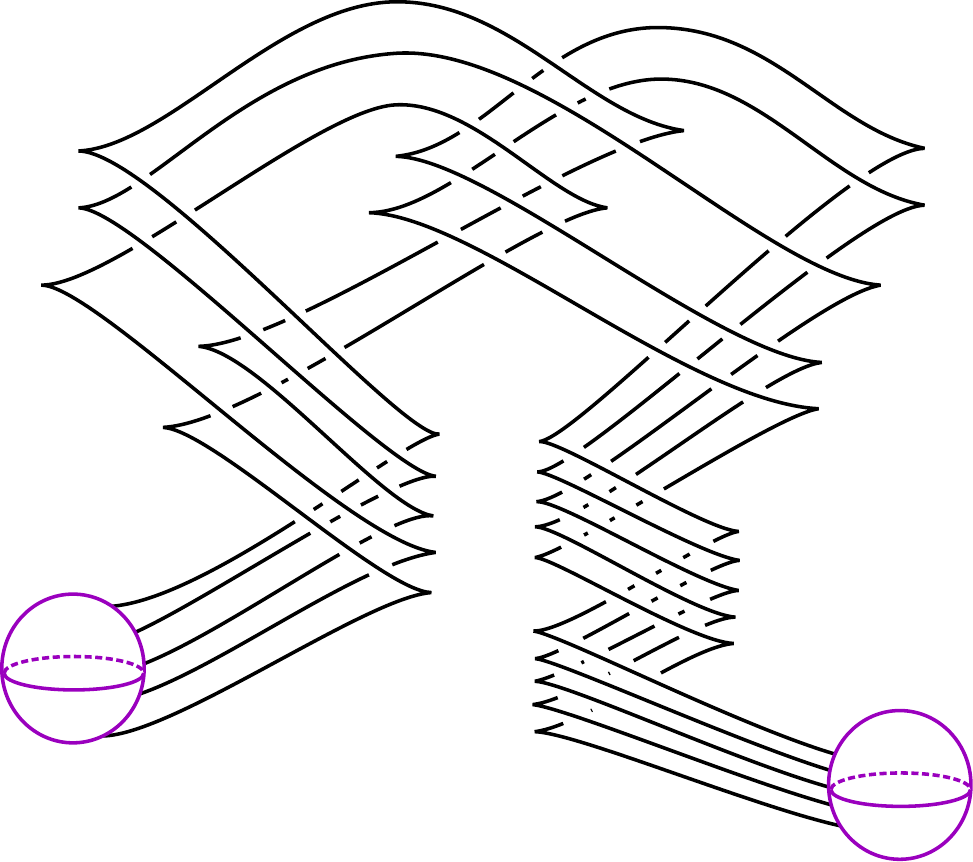}}
                            
                            \put(0.15cm,4.82cm){$-2$}
                        
                        \end{picture} 
                        }
                        \caption{A Legendrian diagram for $C(1,1;-1)$ whose 2-handle has framing $-2$ and Thurston-Bennequin number $-1$. \vspace{-5pt}}
                        \label{20240420-23}
            \end{figure}

    Finally, we slide all five strands of the $2$-handle over the $1$-handle. This decreases both the framing and the Thurston-Bennequin number by $2$, so the immediate resulting front has framing $-2$ and Thurston-Bennequin number $-6$. Moreover, this operation creates a $(-2)$-twist box involving all five strands going over the 1-handle. Absorbing one of these full twists into the five cusp pairs and simplifying the front yields the Legendrian diagram in Figure \ref{20240420-23}. Note that in this diagram, the $(-2)$-twist box is displayed as a $(-\frac{1}{2})$-twist on the left and a $(-\frac{3}{2})$-twist on the right.

    Now the Thurston-Bennequin number of the diagram is $-1$, which is larger than the framing, now equal to $-2$ after the handle slide. Thus, it follows from \cite{Gom98} that $C(1,1;-1)$ admits a Stein structure.

\section{An infinite family of infinite order Stein corks} {\label{section-modified-stein}}

    In this section, we prove that the modified corks $C_m$ all admit Stein structures.

\begin{figure}[H]
                        \centering
                        \resizebox{15.3cm}{!}{
                        \begin{picture}(17cm,5cm)
                            \put(0cm,0.0cm){\includegraphics[width=7cm, height=5cm]{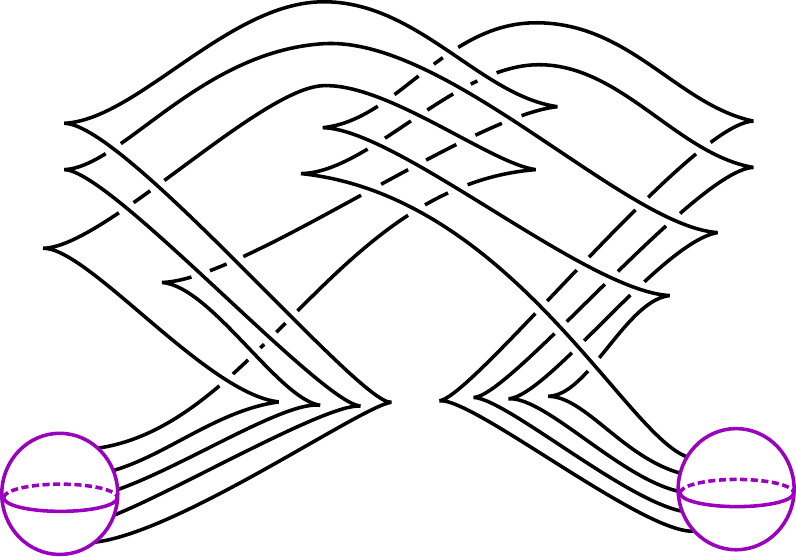}}
                            \put(8cm,0cm){\includegraphics[width=9cm, height=5cm]{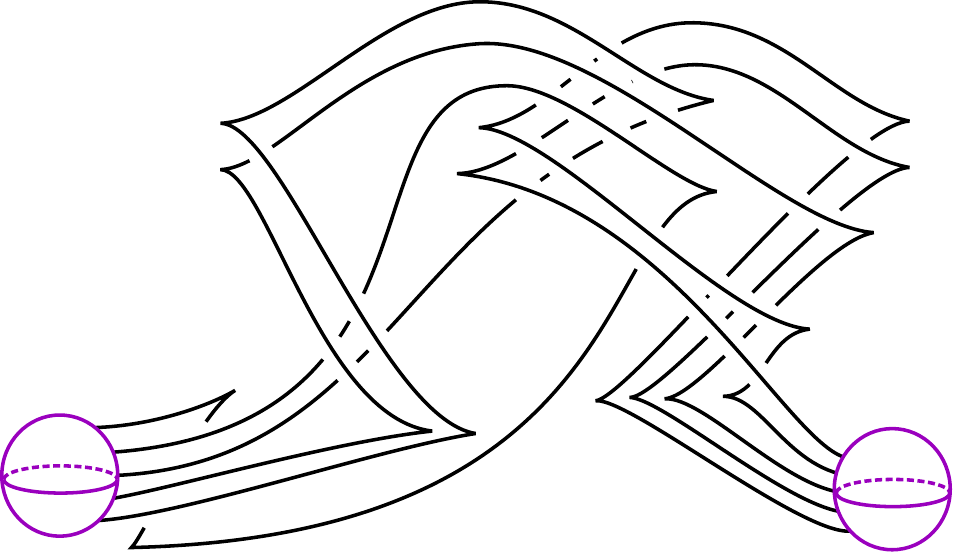}}

                            \put(6.35cm,4.32cm){$0$}
                            \put(16.05cm,4.32cm){$0$}
                        
                        \end{picture}  
                        }
                        \caption{Legendrian diagrams for the cork $C(1,1;-1)$. The left one has Thurston-Bennequin number $-2$, and the right one $-1$.}
                        \label{20240420-25-26}
            \end{figure}

    We start by modifying the handle diagram for $C(1,1;-1)$ given by the second diagram in Figure \ref{20240420-19-20}. A smooth isotopy, shown in the following diagram, gives the right diagram of Figure \ref{20240420-25-26}, in which the black curve has Thurston-Bennequin number $-1$.

    Next, we follow the procedure in Section \ref{nutshell} and take the satellite of the purple $1$-handle (treated as a dotted curve) using Yasui's patterns $P_m$ (cf. Figures \ref{20240604-2} and \ref{pattern-Yasui}).    

\begin{figure}[H]
                        \centering
                        \resizebox{10.8cm}{!}{
                        \begin{picture}(12cm,6cm)
                            \put(0cm,0.5cm){\includegraphics[width=5.5cm, height=4.5cm]{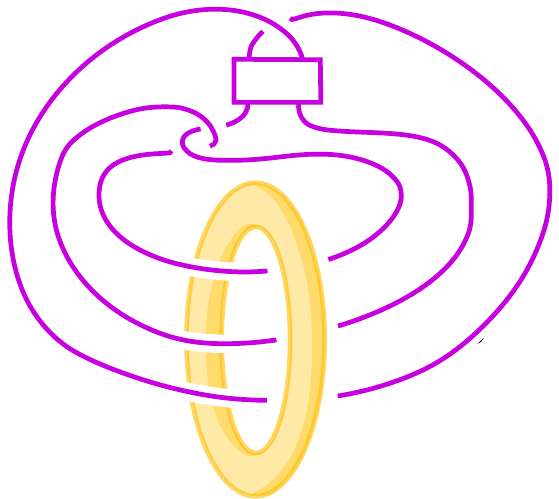}}
                            \put(6.5cm,0cm){\includegraphics[width=5.5cm, height=6cm]{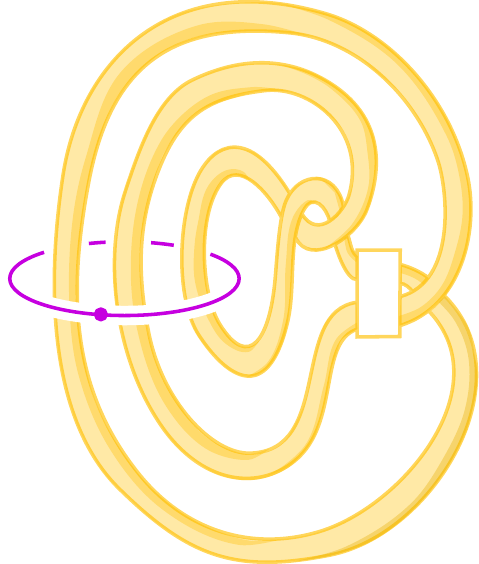}}
                            \put(2.5cm,4.24cm){\scriptsize{$-m$}}
                            \put(10.6cm,2.82cm){\scriptsize{$-m$}}
                        
                        \end{picture}  
                        }
                        \caption{The two ways of applying the pattern $P_m$ are equivalent.}
                        \label{20241030-2}
            \end{figure}

    However, observe that the patterns $P_m$ are symmetric in the following sense. If a link is contained in the union of two Hopf-linked solid tori, then replacing either solid torus component by its satellite with the pattern $P_m$ produces isotopic links. This symmetry is illustrated in Figure \ref{20241030-2}.

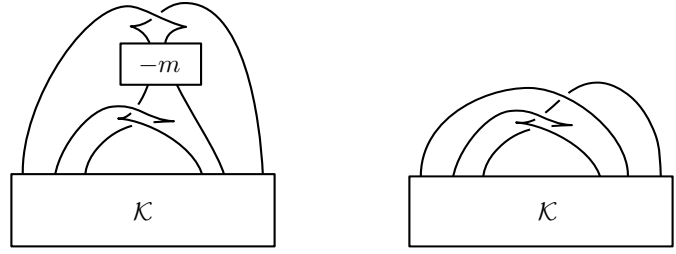
\begin{figure}[H]
                        \centering
                        \resizebox{9cm}{!}{
                        \begin{picture}(10cm,3.7cm)
                            \put(0cm,0.0cm){\includegraphics[width=4cm, height=3.7cm]{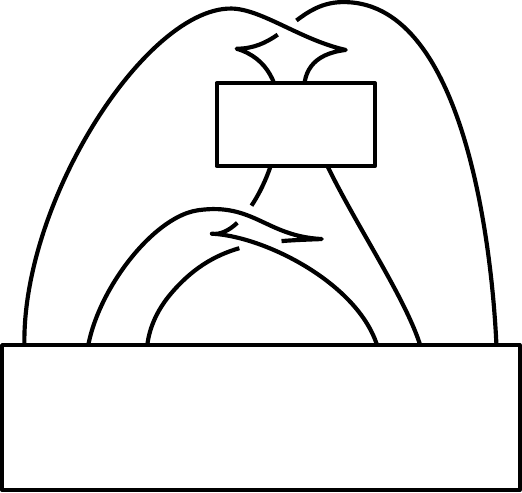}}
                            \put(6cm,0cm){\includegraphics[width=4cm, height=2.5cm]{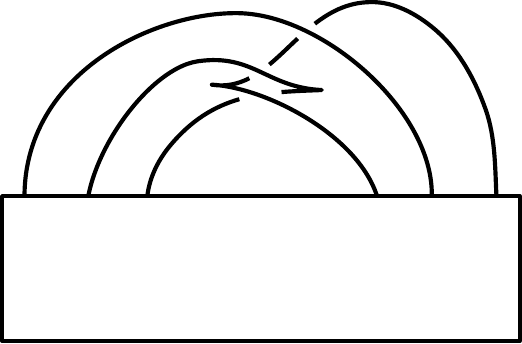}}
                            \put(1.93cm,2.65cm){$-m$}

                            \put(1.85cm,0.42cm){$\mathcal{K}$}
                            \put(7.95cm,0.42cm){$\mathcal{K}$}
                        
                        \end{picture}   
                        }
                        \caption{A Legendrian representative obtained by applying Yasui's $P_m$ to a Legendrian knot $\mathcal{K}$. Left: $m>0$. Right: $m=0$.}
                        \label{20260605-1}
            \end{figure}

\begin{figure}[H]
                        \centering
                        \resizebox{10.8cm}{!}{
                        \begin{picture}(12cm,8cm)
                            \put(0cm,0cm){\includegraphics[width=12cm, height=8cm]{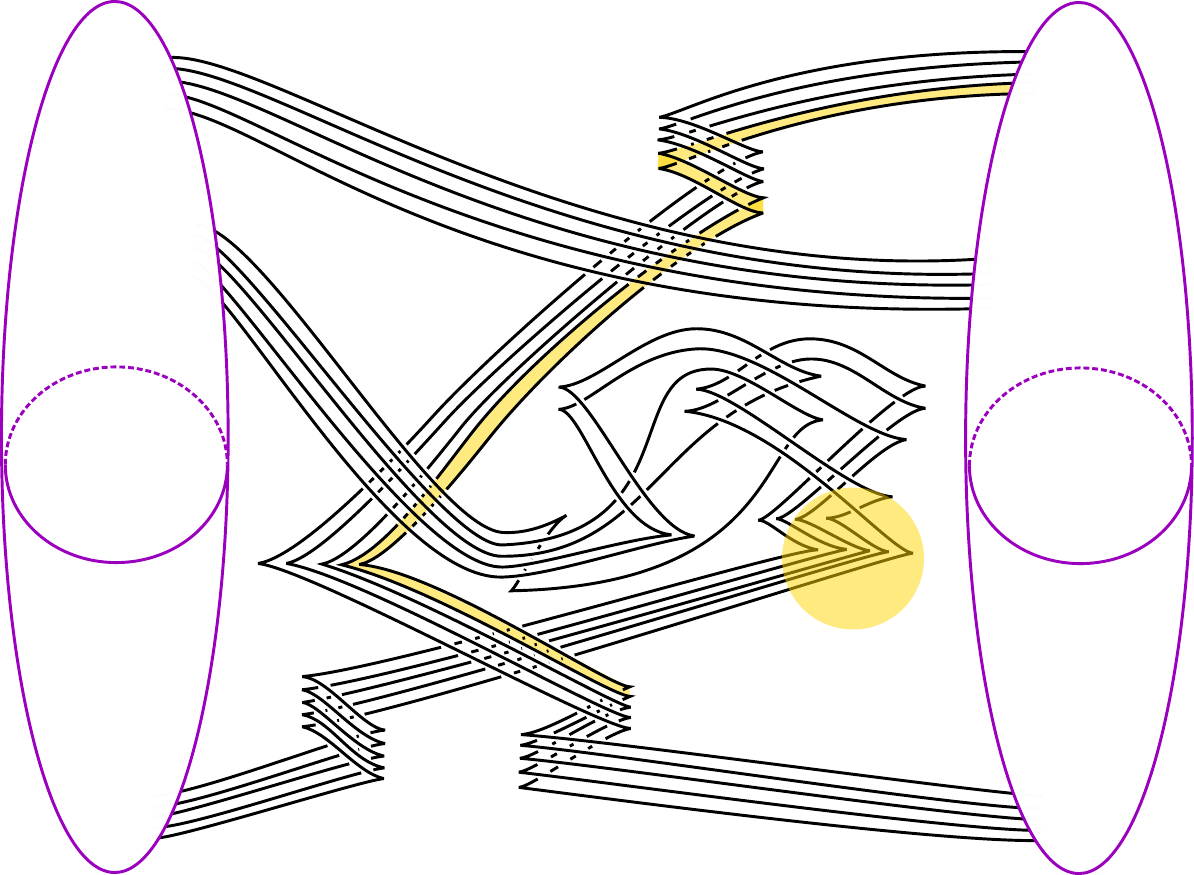}}
                            
                            \put(3.45cm,7.10cm){$0$}
                        
                        \end{picture}    
                        }
                        \caption{A Legendrian diagram for the cork $C_0$, whose 2-handle has framing $0$ and Thurston-Bennequin number $-5$. }
                        \label{20240420-28}
            \end{figure}

    Hence, we can obtain a Legendrian diagram of $C_m$ by taking the Legendrian satellite of the $2$-handle, viewed as a knot in $S^1\times S^2$. Figure \ref{20260605-1} (originally from \cite{Yas17}) shows a Legendrian representative of the pattern $P_m$ applied to a Legendrian knot $\mathcal{K}$. In the case $m=0$, applying this pattern to the black curve in the right-hand diagram of Figure \ref{20240420-25-26} produces the diagram in Figure \ref{20240420-28}, in which the Thurston-Bennequin number of the black curve is $-5$. 

            \begin{figure}[H]
                        \centering
                        \resizebox{10.8cm}{!}{
                        \begin{picture}(12cm,8cm)
                            \put(0cm,0cm){\includegraphics[width=12cm, height=8cm]{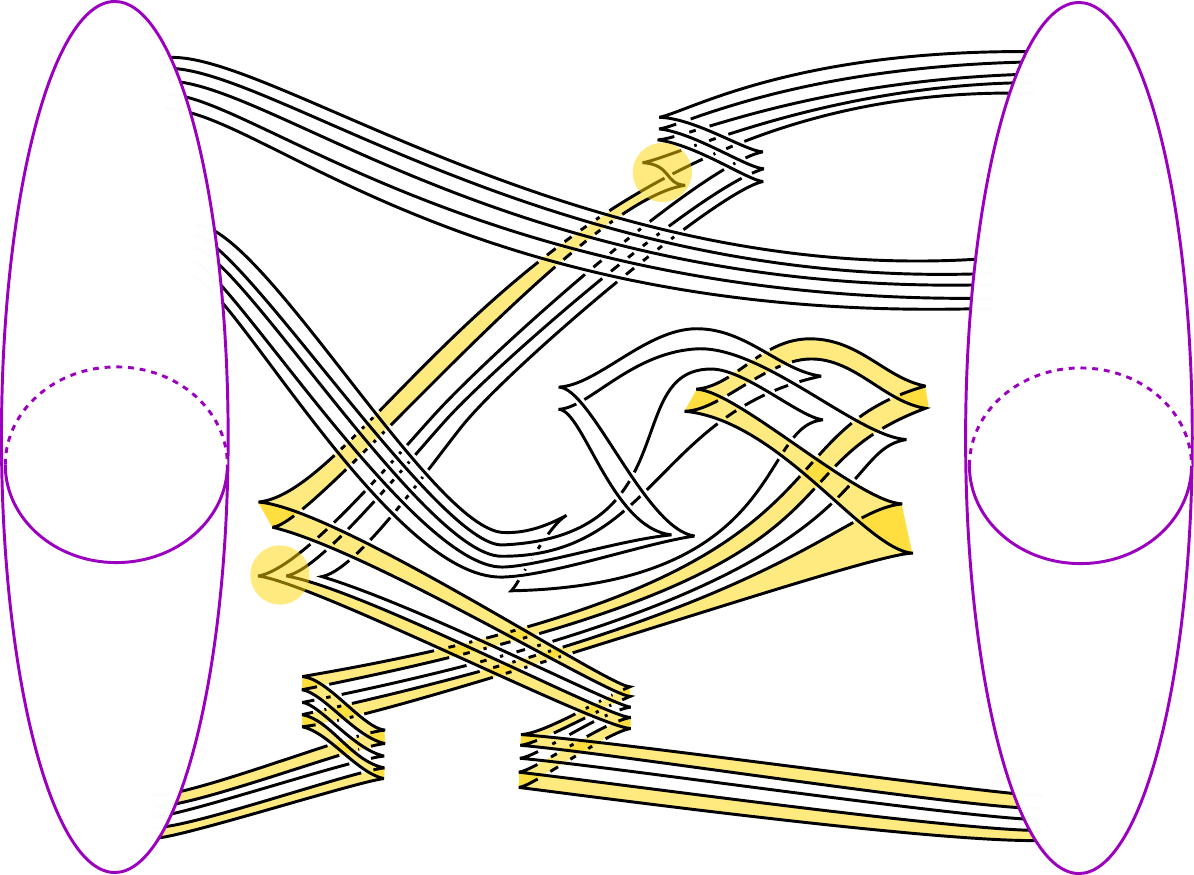}}
                            
                            \put(3.45cm,7.10cm){$0$}
                        
                        \end{picture}    
                        }
                        \caption{A Legendrian diagram for the cork $C_0$, whose 2-handle has framing $0$ and Thurston-Bennequin number $-1$. }
                        \label{20241031-1}
            \end{figure}
            
    Next, we move the two zigzags along the highlighted strip on the left, and straighten the four zigzags in the highlighted region on the right. This gives Figure \ref{20241031-1}, in which the $2$-handle curve has Thurston-Bennequin number $-1$.

\begin{figure}[H]
                        \centering
                        \resizebox{10.8cm}{!}{
                        \begin{picture}(12cm,8cm)
                            \put(0cm,0cm){\includegraphics[width=12cm, height=8cm]{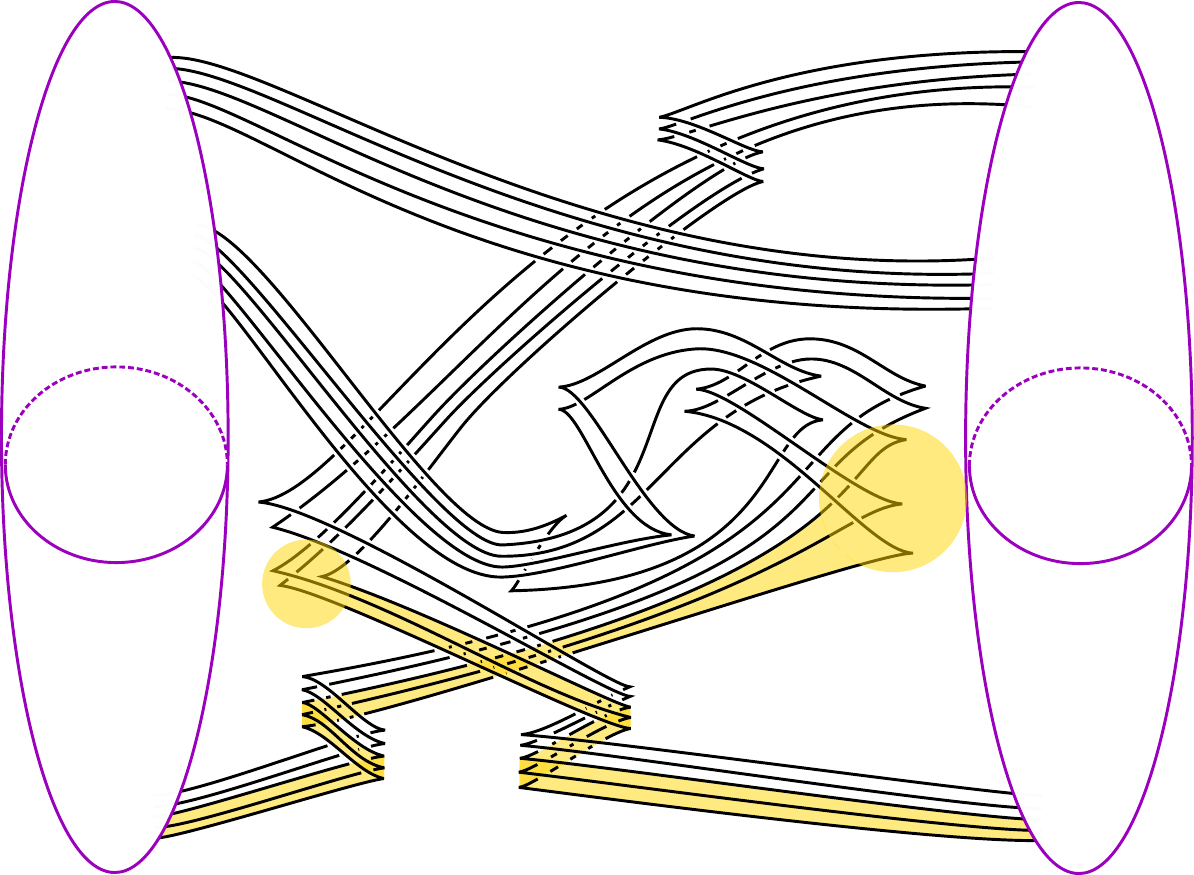}}
                            
                            \put(3.45cm,7.10cm){$0$}
                        
                        \end{picture} 
                        }
                        \caption{A Legendrian diagram for the cork $C_0$, whose 2-handle has framing $0$ and Thurston-Bennequin number $0$. }
                        \label{20241222-1}
            \end{figure}

    Then we move the highlighted zigzag at the top along the strip. This gives Figure \ref{20241222-1}. This step preserves the writhe while decreasing the number of left cusps by $1$; hence the Thurston-Bennequin number is increased by $1$ and is now $0$.

\begin{figure}[H]
                        \centering
                        \resizebox{10.8cm}{!}{
                        \begin{picture}(12cm,8cm)
                            \put(0cm,0cm){\includegraphics[width=12cm, height=8cm]{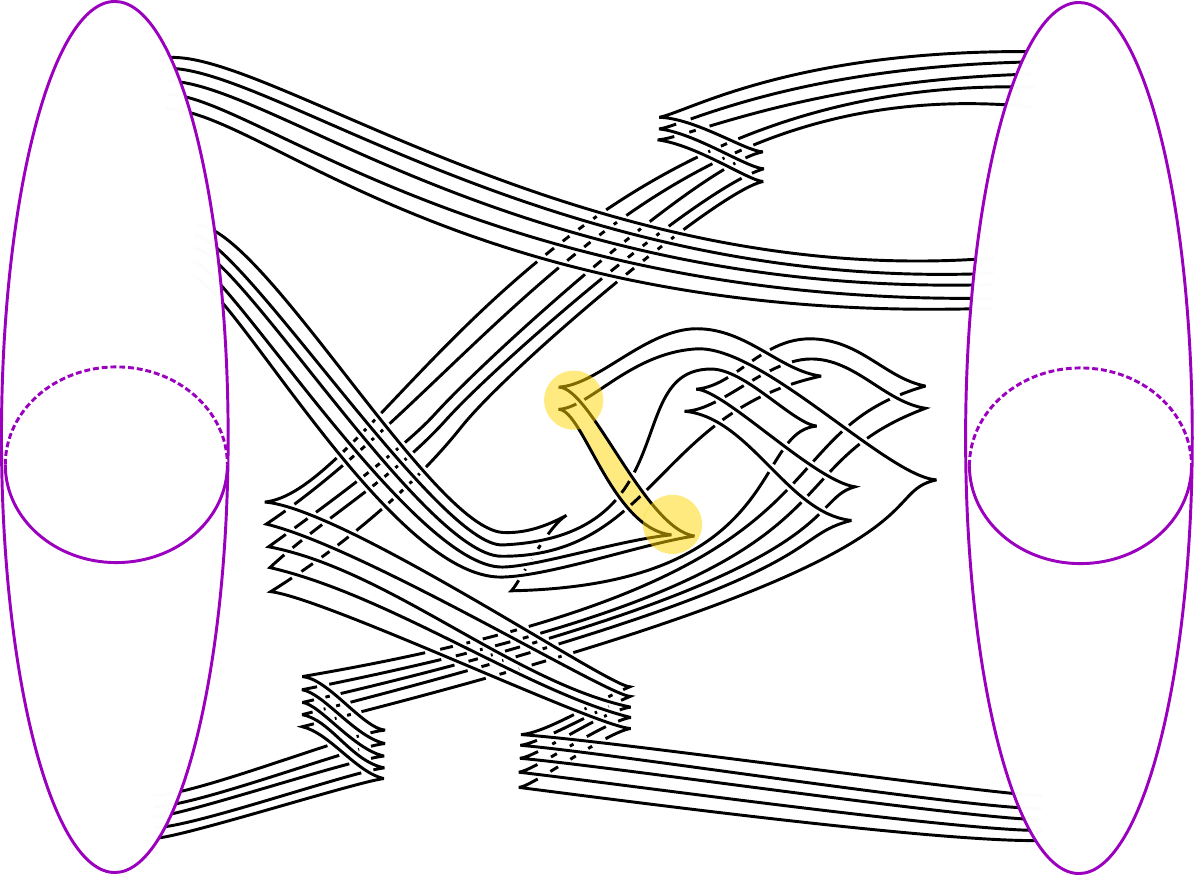}}
                            
                            \put(3.45cm,7.10cm){$0$}
                        
                        \end{picture} 
                        }
                        \caption{A Legendrian diagram for the cork $C_0$, whose 2-handle has framing $0$ and Thurston-Bennequin number $0$. }
                        \label{20241222-2}
            \end{figure}

\begin{figure}[H]
                        \centering
                        \resizebox{10.8cm}{!}{
                        \begin{picture}(12cm,8cm)
                            \put(0cm,0cm){\includegraphics[width=12cm, height=8cm]{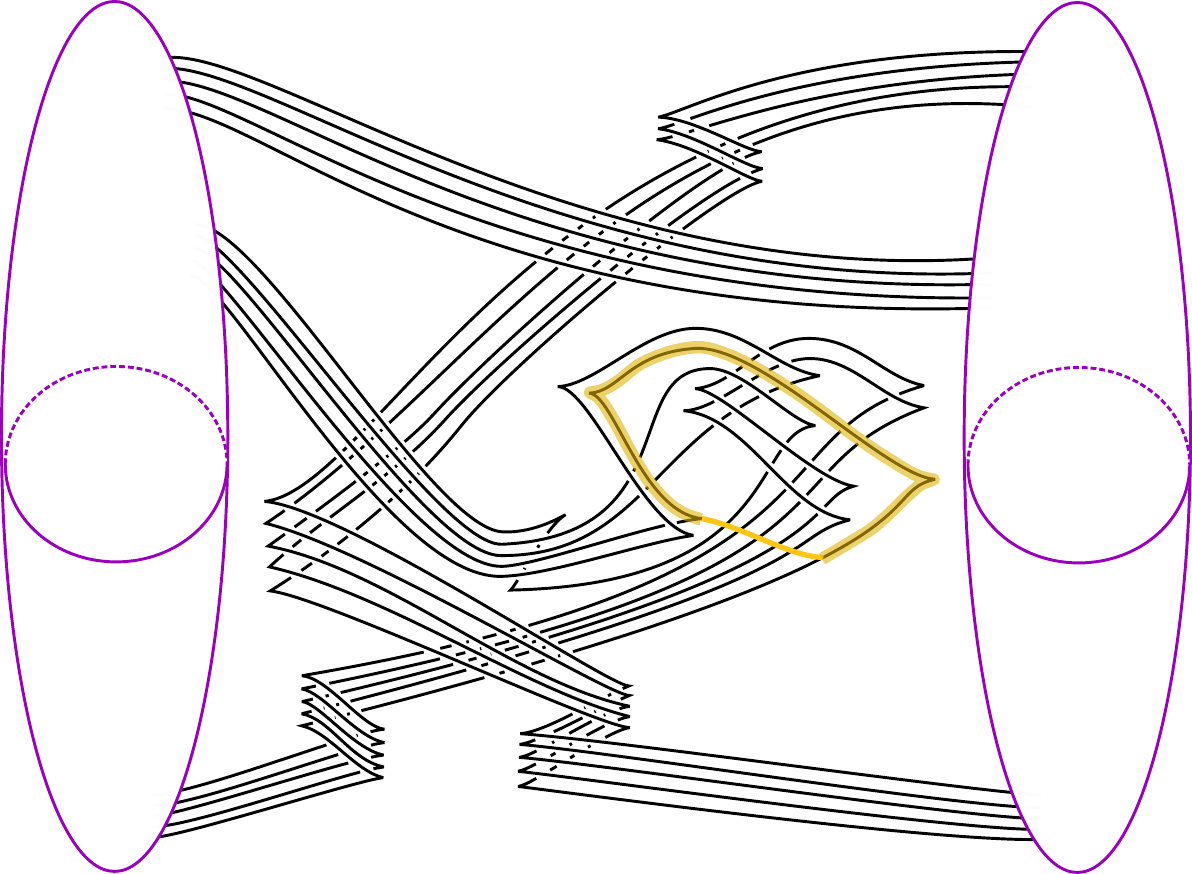}}
                            
                            \put(3.45cm,7.10cm){$0$}
                        
                        \end{picture}   
                        }
                        \caption{A Legendrian diagram for the cork $C_0$, whose 2-handle has framing $0$ and Thurston-Bennequin number $0$. }
                        \label{20241222-3}
            \end{figure}

        Next, in Figure \ref{20241222-1}, we again move the highlighted zigzag on the left along the strip. This gives Figure \ref{20241222-2}. A similar operation gives Figure \ref{20241222-3}. Since both steps can be achieved via Legendrian isotopies, the Thurston-Bennequin number remains unchanged.

\begin{figure}[H]
                        \centering
                        \resizebox{10.8cm}{!}{
                        \begin{picture}(12cm,8cm)
                            \put(0cm,0cm){\includegraphics[width=12cm, height=8cm]{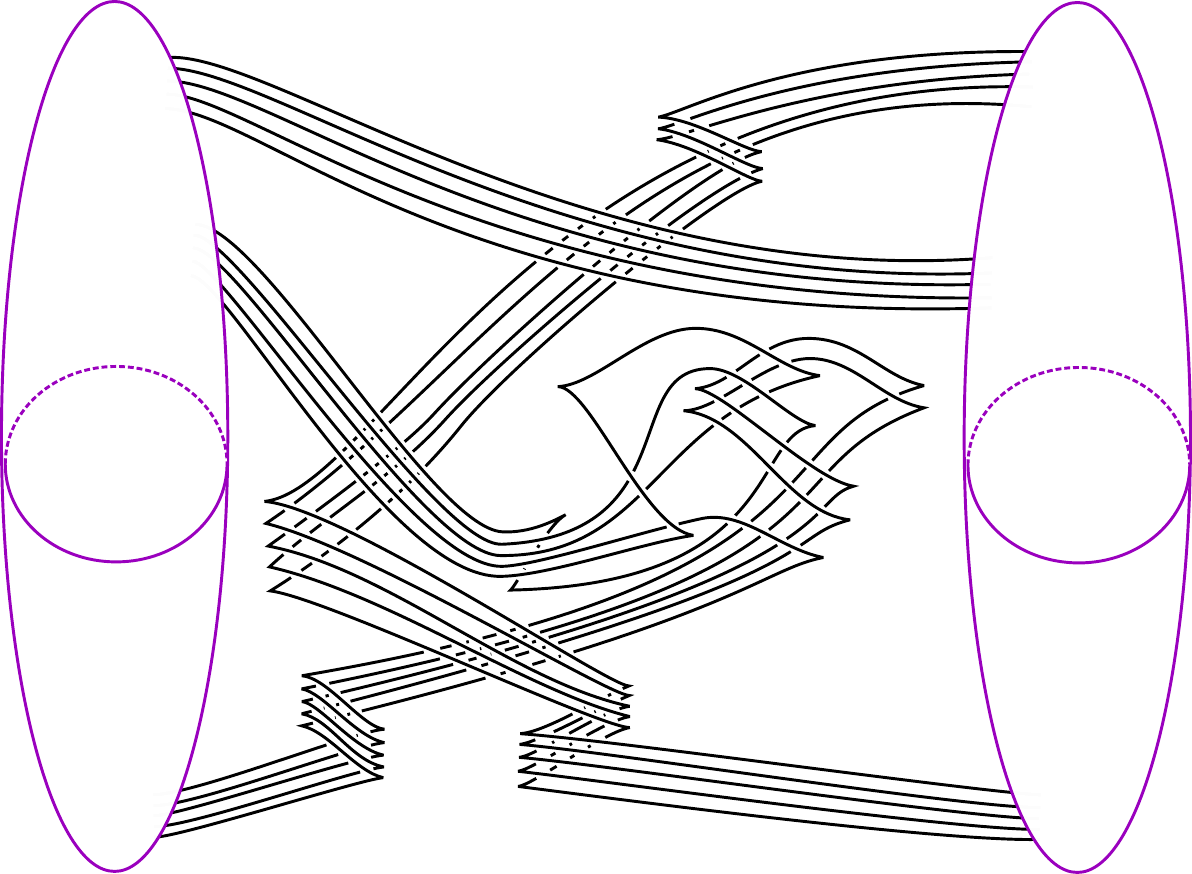}}
                            
                            \put(3.45cm,7.10cm){$0$}
                        
                        \end{picture}    
                        }
                        \caption{A Legendrian diagram for the cork $C_0$, whose 2-handle has framing $0$ and Thurston-Bennequin number $+1$. }
                        \label{20241222-4}
            \end{figure}

        Finally, we pull the highlighted strand down to the orange replacement curve, as indicated in Figure \ref{20241222-3}. This gives Figure \ref{20241222-4}. Observe that the Thurston-Bennequin number increases by $1$, and is therefore now larger than the framing of the $2$-handle (which is $0$).

\begin{figure}[H]
                        \centering
                        \resizebox{10.8cm}{!}{
                        \begin{picture}(12cm,9cm)
                            \put(0cm,0cm){\includegraphics[width=12cm, height=9cm]{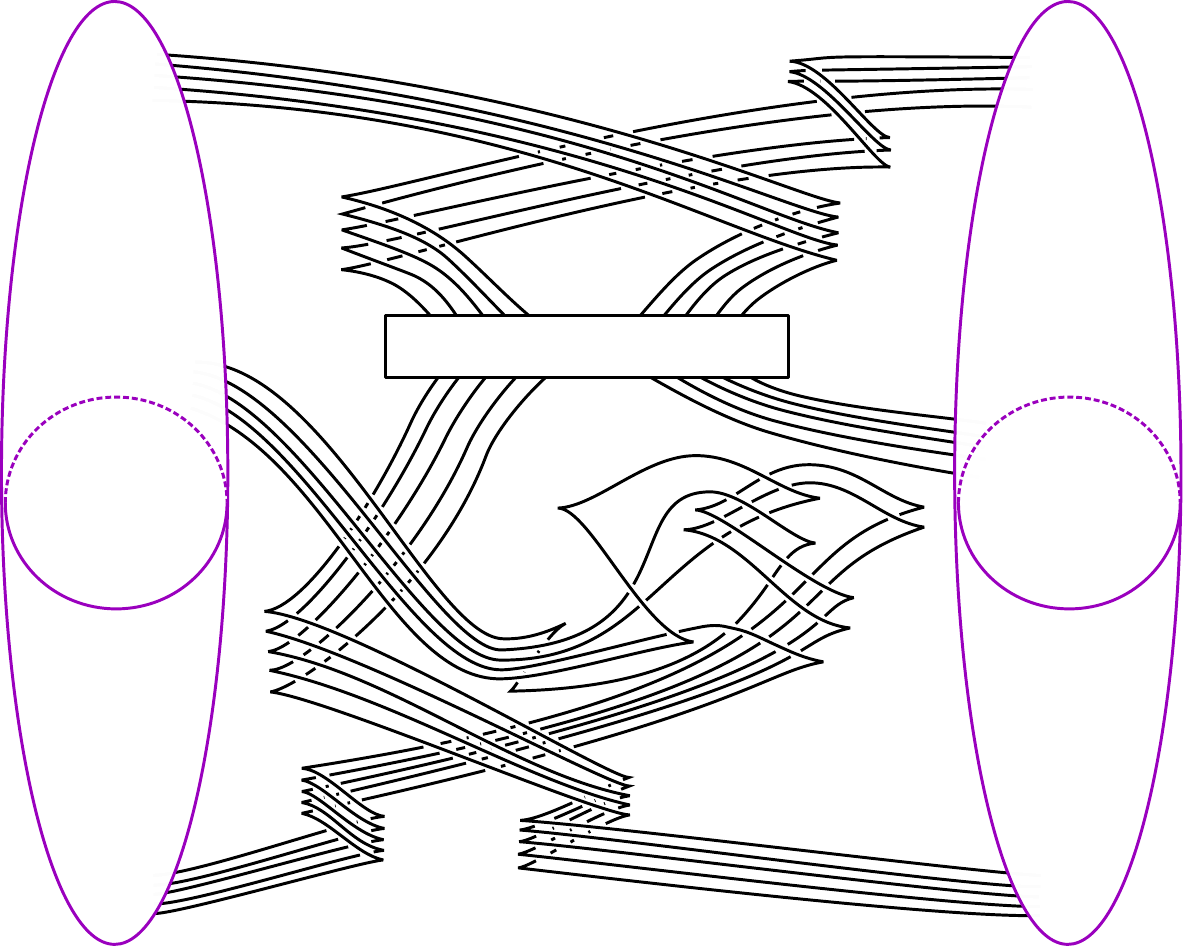}}
                            
                            \put(5.62cm,5.57cm){$-m$}
                            \put(3.45cm,8.40cm){$0$}
                        
                        \end{picture}  
                        }
                        \caption{A Legendrian diagram for the cork $C_m$, $m>0$, whose 2-handle has framing $0$ and Thurston-Bennequin number $+1$. }
                        \label{20241222-5}
            \end{figure}

        This shows that the manifold $C_0$ admits a Stein structure. To prove that each $C_m$ is a Stein domain for $m>0$, we modify the above diagram to the one shown in Figure \ref{20241222-5}. Observe that the part of the diagram outside the $(-m)$-twist box contributes $+1$ to the Thurston-Bennequin number, as before.

\begin{figure}[H]
                        \centering
                        \resizebox{6.4cm}{!}{
                        \begin{picture}(8cm,5cm)
                            \put(0cm,0cm){\includegraphics[width=3cm, height=5cm]{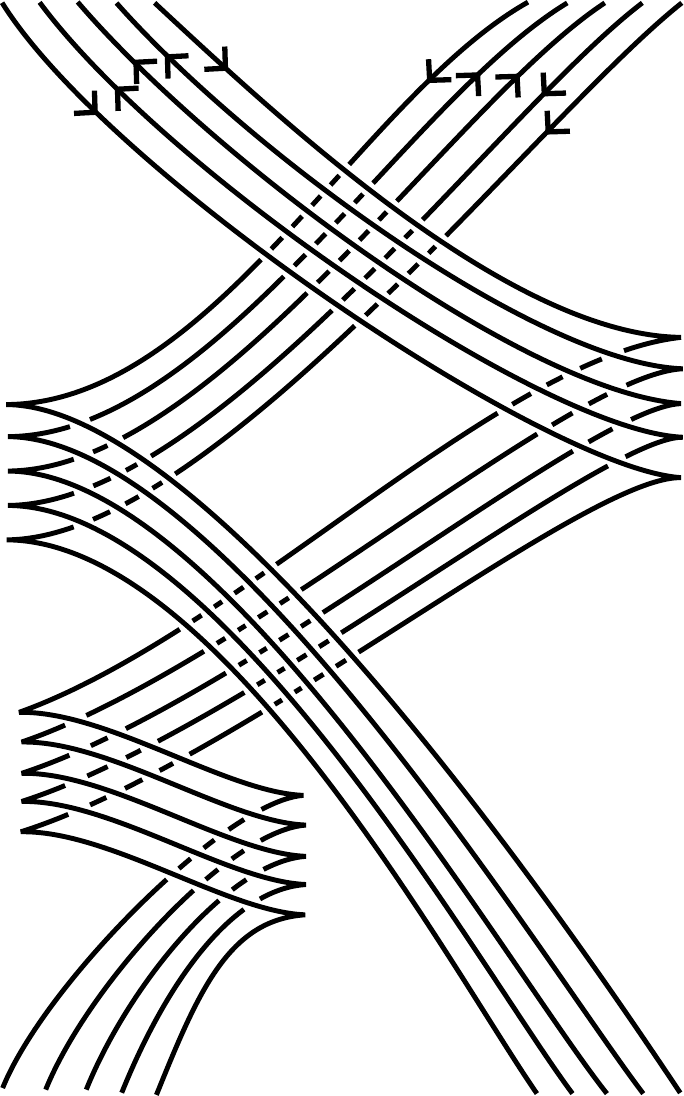}}

                            \put(5cm,0cm){\includegraphics[width=3cm, height=5cm]{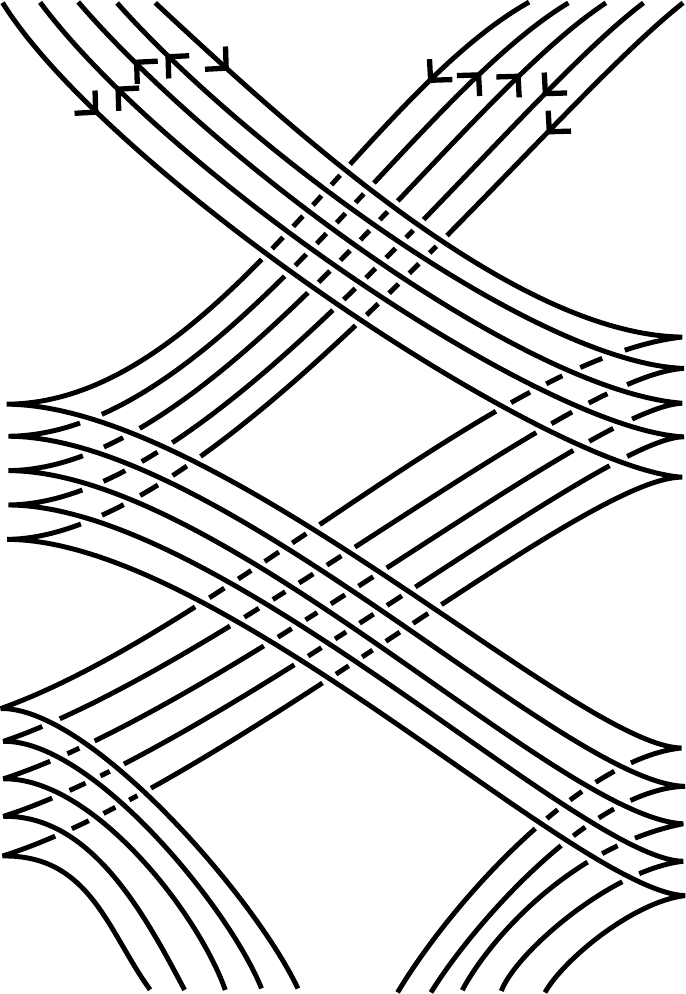}}
                        
                        \end{picture} 
                        }
                        \caption{Two local models for a negative full twist.}
                        \label{20241223-1}
            \end{figure}

        Inside the $(-m)$-twist box, each negative twist corresponds to one of the local diagrams depicted in Figure \ref{20241223-1}: the bottom twist corresponds to the diagram on the left, while the remaining twists correspond to the diagram on the right. In both local diagrams, the contribution to the Thurston-Bennequin number is $0$, which implies that the $(-m)$-twist box does not change the Thurston-Bennequin number.

        This completes the proof of the main theorem. \hfill\raisebox{-0.25ex}{\qedsymbol}

\medskip

\bibliographystyle{alpha} 
\bibliography{bib}

\end{spacing}
\end{document}